\documentclass[11pt,a4paper,reqno]{article}
\usepackage[a4paper,top=25mm,bottom=25mm,left=26mm,right=26mm]{geometry}
\usepackage{amsmath,amssymb,amsthm,mathtools,bm}
\usepackage[T1]{fontenc}
\usepackage[hidelinks]{hyperref}

\numberwithin{equation}{section}
\allowdisplaybreaks[1]
\newtheorem{theorem}{Theorem}[section]
\newtheorem{proposition}[theorem]{Proposition}
\newtheorem{lemma}[theorem]{Lemma}
\newtheorem{corollary}[theorem]{Corollary}

\newtheoremstyle{boldremark}
  {\topsep}
  {\topsep}
  {\normalfont}
  {}
  {\bfseries}
  {.}
  {0.5em}
  {}
\theoremstyle{boldremark}
\newtheorem{remark}[theorem]{Remark}

\newcommand{\R}{\mathbb R}
\newcommand{\dd}{\,dx}
\newcommand{\ip}[2]{\langle #1,#2\rangle}
\newcommand{\norm}[1]{\lVert #1\rVert}
\newcommand{\Hm}{\dot H^{-1}}
\newcommand{\Hp}{\dot H^1}
\newcommand{\bu}{\bar u_1}
\newcommand{\br}{\bar\rho}
\newcommand{\bt}{\bar\theta}
\newcommand{\rr}{\widetilde\rho}
\newcommand{\EE}{\widetilde E}
\newcommand{\mm}{\widetilde m}
\newcommand{\curl}{\operatorname{curl}}
\newcommand{\diver}{\operatorname{div}}
\newcommand{\diag}{\operatorname{diag}}
\newcommand{\cA}{\mathcal A}
\newcommand{\cL}{\mathcal L}
\newcommand{\cN}{\mathcal N}
\newcommand{\cE}{\mathcal E}

\newcommand{\ETan}{\widetilde{\mathcal E}_1}

\hypersetup{
  pdftitle={Nonlinear Stability of Planar Shock Waves for the Navier-Stokes-Fourier Equations in R\string^3},
  pdfauthor={Feimin Huang and Lingda Xu}
}

\begin{document}

\begin{center}
{\Large\bfseries Nonlinear Stability of Planar Shock Waves for the\\[4pt]
Navier--Stokes--Fourier Equations in $\R^3$}\\[14pt]
{\large Feimin Huang$^{1}$ \quad and \quad Lingda Xu$^{2}$}\\[8pt]
{\small
$^{1}$Academy of Mathematics and Systems Science, Chinese Academy of Sciences, Beijing 100190, China.\\[2pt]
$^{2}$Hetao Institute of Mathematics and Interdisciplinary Sciences, Shenzhen 518048, China.
}
\end{center}
\begingroup
\renewcommand\thefootnote{}\footnote{E-mail addresses: \texttt{fhuang@amt.ac.cn} (F. Huang), \texttt{xuldmath@gmail.com} (L. Xu).}%
\addtocounter{footnote}{-1}%
\endgroup

\begin{abstract}
In this paper, we investigate the global nonlinear stability and large-time decay rates of small-amplitude planar viscous shock waves for the three-dimensional (3-D) full compressible Navier--Stokes--Fourier equations in the whole space $\R^3$. {We prove that if the initial perturbation is sufficiently small in $H^3(\R^3)\cap\Hm(\R^3)$, then the Cauchy problem admits a unique global strong solution that converges to the planar shock profile without any phase shift. In the moving shock frame, both the conservative and primitive perturbations decay at the rate $(1+t)^{-1/4}$ in $H^3(\R^3)$. Their time and transverse spatial derivatives decay at the rate $(1+t)^{-1/2}$, measured in $\Hm(\R^3)$ for the conservative perturbation and in $H^1(\R^3)$ for the primitive perturbation.} No zero-mass condition or spectral stability assumption is imposed. On the whole space $\R^3$, the main analytical difficulties stem from the continuous transverse low frequencies, the opposite signs of shock compression in $L^2(\R^3)$ and unweighted $\Hm(\R^3)$, and the absence of density dissipation in the continuity equation coupled with temperature variations. To overcome these difficulties, we introduce a density-weighted $\Hm(\R^3)$ multiplier built from the inverse entropy Hessian and a modified curl operator. Through entropy symmetrization and the curl projection, all second-order acoustic derivatives of the potential cancel out and the transverse acoustic modes are controlled, so that the normal Euler drift is governed by the compressive derivative of the characteristic speed. Combined with an effective momentum variable and a momentum--density cross estimate, this negative-order energy absorbs the unfavorable compression term in the $L^2$ energy and closes the $H^3(\R^3)$ a priori estimates without loss of derivatives.
\end{abstract}



\section{Introduction}\label{sec:intro}

\leavevmode{Consider the three-dimensional full compressible Navier--Stokes--Fourier equations in Eulerian coordinates $(t,y)\in(0,\infty)\times\R^3$:}
\begingroup\begin{equation}\label{eq:nsf}
\left\{
\begin{aligned}
&\rho_t+\diver(\rho u)=0,\\
&(\rho u)_t+\diver(\rho u\otimes u)+\nabla p=\cL u,\\
&E_t+\diver((E+p)u)=\diver\{S(u)u+\kappa\nabla\theta\},
\end{aligned}
\right.\qquad y\in\R^3,\ t>0,
\end{equation}
\endgroup
where $\rho>0$, $u=(u_1,u_2,u_3)^T\in\R^3$, and $\theta>0$ denote the fluid density, velocity, and absolute temperature, respectively. For an ideal polytropic gas, the pressure $p$, total energy $E$, viscous stress tensor $S(u)$, and Lam\'e operator $\cL$ are given by
\begin{gather*}
p=\rho\theta,\qquad E=\rho\left(c_v\theta+\frac12|u|^2\right),\qquad \gamma=1+c_v^{-1}>1,\\
S(u)=\mu(\nabla u+\nabla u^T)+\lambda(\diver u)I,
\qquad\cL u=\mu\Delta u+(\mu+\lambda)\nabla\diver u.
\end{gather*}
Throughout the paper, the viscosity coefficients $\mu>0$ and $\lambda$ (with $2\mu+3\lambda\ge0$), the specific heat $c_v>0$, and the heat conductivity $\kappa>0$ are fixed constants, and we set $\nu=2\mu+\lambda>0$. {Our goal is to establish the global nonlinear stability and explicit large-time decay rates of a small-amplitude planar viscous shock wave $\bar U(y_1-st)$ in $\R^3$, where $U=(\rho,\rho u,E)^T$ and $\bar u=\bu e_1$ with $\bu'<0$.}

The stability of elementary waves---viscous shock waves, rarefaction waves, and contact discontinuities---for one-dimensional viscous conservation laws has been extensively studied since the 1980s. For 1-D viscous shock profiles, Matsumura--Nishihara~\cite{mn} and Goodman~\cite{goodman} independently introduced the anti-derivative method and exploited the compressive monotonicity $\bu'<0$ to prove nonlinear stability under the zero-mass condition; see also Kawashima--Matsumura~\cite{km}. {When the zero-mass condition is removed, perturbations carry non-trivial masses along trans\-versal characteristic families; Liu~\cite{liu85,liu97} overcame this difficulty by introducing coupled diffusion waves, and Szepessy--Xin~\cite{sx} established nonlinear stability via spatially weighted energy estimates. Pointwise Green's function bounds and detailed wave-interaction descriptions were subsequently developed by Liu--Zeng for general hyperbolic--parabolic systems~\cite{lz97} and physical gas dynamics~\cite{lz15}, while spectral and nonlinear orbital stability of small-amplitude shocks were obtained by Humpherys--Zumbrun~\cite{hz} and Mascia--Zumbrun~\cite{mz}.} For 1-D rarefaction waves, contact discontinuities, and composite Riemann patterns, we refer to~\cite{mn86,mn92,xin93,hmx,hxy,hlm,lwx}; {in particular, based on the method of $a$-contraction with shifts~\cite{kv17,kv21}, Kang--Vasseur--Wang~\cite{kvw} recently established the time-asymptotic stability of generic Riemann solutions consisting of a shock, a contact discontinuity, and a rarefaction wave for the 1-D compressible Navier--Stokes--Fourier system.}

\leavevmode{In multiple space dimensions, the global well-posedness and asymptotic decay of perturbations around constant equilibrium states were pioneered by Matsumura--Nishida~\cite{matsunishi,mn83} and further developed in the framework of dissipative hyperbolic--parabolic systems~\cite{kawashima,ponce}, multidimensional diffusion waves~\cite{hoffzumbrun}, scaling-critical Besov spaces~\cite{danchin}, and negative Sobolev interpolation~\cite{gw}. For multidimensional planar rarefaction waves, nonlinear stability was established by Xin~\cite{xin90} for scalar viscous conservation laws, Hokari--Matsumura~\cite{hm} for the 2-D Euler system with artificial viscosity, and Li--Wang~\cite{lw18} and Li--Wang--Wang~\cite{lww20} for the isentropic and full compressible Navier--Stokes equations on periodic cylinders. Related stability results on $\R\times\mathbb T^2$ were recently obtained by Duan--Huang--Li--Xu for planar entropy waves in the compressible Navier--Stokes system~\cite{dhlx25} and the Landau equation~\cite{dhlx26}.}

\leavevmode{For multidimensional planar viscous shocks, Zumbrun~\cite{zumbrun} developed a general spectral framework relating linearized and nonlinear stability to an Evans-function condition; under this spectral hypothesis, Nguyen~\cite[Theorem~1.8]{nguyen} proved $L^p$ decay ($2\le p\le\infty$) of Lax shocks in dimensions $d\ge3$ for sufficiently regular $L^1\cap H^s$ perturbations. Spectral stability has been verified analytically by Freist\"uhler--Szmolyan~\cite{fs} for small shocks with artificial Laplacian viscosity, and numerically combined with high-frequency analytical bounds by Humpherys--Lyng--Zumbrun~\cite{hlz} for physical Navier--Stokes shocks.}

Recently, Wang--Wang~\cite{ww} proved the nonlinear stability of small planar viscous shocks for the 3-D isentropic Navier--Stokes system on the periodic domain $\R\times\mathbb T^2$ under $H^2$ perturbations up to a time-dependent shift, without imposing a spectral hypothesis or a zero-mass condition. Whether planar viscous shocks are stable for the full non-isentropic Navier--Stokes--Fourier system \eqref{eq:nsf} in the whole space $\R^3$, and at what rate perturbations decay in time, is the problem we address in this paper.

\leavevmode{Let $K\subset(0,\infty)\times\R^3\times(0,\infty)$ be a fixed compact set in the physical state space, and consider two endstates $(\rho_\pm,u_\pm,\theta_\pm)\in K$ with $u_\pm=u_{1,\pm}e_1$ connected by a right-moving (3-family) Lax shock. Their conservative states $U_\pm=(\rho_\pm,\rho_\pm u_\pm,E_\pm)^T$ and the shock speed $s$ satisfy the Rankine--Hugoniot relation and Lax entropy condition}
\[
s[U]=[F_1(U)],\qquad
\lambda_+(U_+)<s<\lambda_+(U_-),\qquad
\lambda_+(U)=u_1+\sqrt{\gamma\theta},
\]
where $[a]=a_+-a_-$ denotes the jump across the shock. The Euler fluxes are given by
\begin{equation}\label{eq:Eulerflux}
F_j(U)=\bigl(\rho u_j,\ \rho uu_j+pe_j,\ (E+p)u_j\bigr)^T,\qquad
p=c_v^{-1}\left(E-\frac{|m|^2}{2\rho}\right),\quad m=\rho u.
\end{equation}
Let $\delta=\rho_--\rho_+>0$ be the shock strength. {Up to a constant Galilean translation along $e_1$, we normalize the midpoint reference state $(\rho_*,0,\theta_*)$ by}
\begin{equation}\label{eq:reference}
\rho_*=\frac{\rho_-+\rho_+}{2},\qquad
\theta_*=\frac{\theta_-+\theta_+}{2},\qquad
\frac{u_{1,-}+u_{1,+}}{2}=0.
\end{equation}
\leavevmode{Let $\bar U(y_1-st)$ be the unique small-amplitude planar shock profile connecting $U_-$ to $U_+$ with $\bar u=\bu e_1$, $\bu'<0$, and $\bu(0)=0$, where the prime denotes differentiation with respect to the traveling coordinate (see Lemma~\ref{lem:profile}).}

\leavevmode{Passing to the moving frame $x=y-ste_1$, the planar shock profile becomes stationary, $\bar U=\bar U(x_1)$. With a slight abuse of notation, we still write $U(t,x)=(\rho,\rho u,E)^T(t,x)$ for the solution in the $(t,x)$-coordinates, so that $\partial_t$ at fixed $x$ corresponds to $\partial_t+s\partial_{y_1}$ in the original Eulerian coordinates $(t,y)$.} We denote the conservative perturbation by $f$ and the primitive perturbation by $Z=(\phi,\psi,\zeta)^T$:
\begin{equation}\label{eq:variables}
\begin{gathered}
f=(\rr,\mm,\EE)^T=U-\bar U(x_1),\qquad Z=(\phi,\psi,\zeta)^T,\\
\phi=\log\left(\frac{\rho}{\br}\right),\qquad \psi=u-\bar u,\qquad \zeta=\theta-\bt.
\end{gathered}
\end{equation}
{In the moving coordinates $(t,x)$, the material derivative and time derivative are written as}
\begingroup\[
D_t=\partial_t+(u-se_1)\cdot\nabla,\qquad
f_t=\partial_tf,\qquad Z_t=\partial_tZ.
\]
\endgroup
Now we state our main result.

\begin{theorem}\label{thm:main}
There exist positive constants $\delta_0,\varepsilon_0,C>0$, depending only on $(\mu,\lambda,\kappa,c_v)$ and the compact set $K$, such that if $0<\delta\le\delta_0$ and the initial perturbation $f_0$ satisfies
\[
f_0\in H^3(\R^3)\cap\Hm(\R^3),\qquad
\norm{f_0}_{H^3}+\norm{f_0}_{\Hm}\le\varepsilon\le\varepsilon_0,
\]
{then the Cauchy problem \eqref{eq:nsf} with initial data $U(0,y)=\bar U(y_1)+f_0(y)$ admits a unique global strong solution on $[0,\infty)\times\R^3$ with $\rho$ and $\theta$ uniformly bounded above and away from zero.} In the moving coordinates,
\begin{gather*}
Z,f\in C([0,\infty);H^3(\R^3)),\qquad f\in C([0,\infty);\Hm(\R^3)),\\
\nabla(\psi,\zeta)\in L^2(0,\infty;H^3(\R^3)),
\end{gather*}
and it holds that
\begin{equation}\label{eq:global}
\sup_{t\ge0}\left\{\norm{f(t)}_{\Hm}^2+\norm{Z(t)}_{H^3}^2\right\}
+\int_0^\infty\left\{\norm{Z(t)}_{H^3}^2+\norm{\nabla(\psi,\zeta)(t)}_{H^3}^2\right\}\,dt
\le C\varepsilon^2,
\end{equation}
together with the time-decay estimates
\begin{align}
\norm{f(t)}_{H^3}+\norm{Z(t)}_{H^3}&\le C\varepsilon(1+t)^{-1/4},\label{eq:rate}\\
\norm{f_t(t)}_{\Hm}+\norm{Z_t(t)}_{H^1}&\le C\varepsilon(1+t)^{-1/2},
\label{eq:trate}
\end{align}
for all $t\ge0$, where the time derivative in \eqref{eq:trate} is taken in the moving frame.
\end{theorem}

\begin{remark}\label{rem:data}
For any $v\in L^1(\R^3)\cap L^2(\R^3)$, splitting the Fourier integral at $|\xi|=1$ in $\R^3$ gives
\[
\norm v_{\Hm}^2\le C\norm v_1^2\int_0^1 dr+\norm v_2^2
\le C\bigl(\norm v_1^2+\norm v_2^2\bigr).
\]
Hence Theorem~\ref{thm:main} applies to any small perturbation $f_0\in H^3(\R^3)\cap L^1(\R^3)$ without requiring the zero-mass condition. In particular, for any localized front modulation $\alpha\in C_c^\infty(\R^2)$, the perturbed shock front satisfies
\[
\bar U(x_1-\alpha(x_\perp))-\bar U(x_1)\in H^3(\R^3)\cap L^1(\R^3),
\qquad
\norm{\bar U(\cdot-\alpha)-\bar U}_1
\le\norm{\bar U'}_{L^1(\R)}\norm\alpha_{L^1(\R^2)},
\]
{with non-vanishing normal mass $\int_{\R}\{\bar U(x_1-\alpha(x_\perp))-\bar U(x_1)\}\,dx_1=-\alpha(x_\perp)[U]$. Writing $\bar U(x_1-\eta\alpha(x_\perp))-\bar U(x_1)=-\eta\alpha(x_\perp)\int_0^1\bar U'(x_1-\tau\eta\alpha(x_\perp))\,d\tau$ and using the exponential decay of $\bar U'$ and its derivatives from Lemma~\ref{lem:profile}, we see that the perturbation norm in $H^3(\R^3)\cap L^1(\R^3)$ is $O(\eta)$ as $\eta\to0$, so Theorem~\ref{thm:main} applies for all sufficiently small $\eta$. A nonzero constant translation of the entire planar profile does not belong to the present $L^2(\R^3)$ perturbation class. The no-shift conclusion is therefore understood within this localized perturbation class, which nevertheless includes sufficiently small localized front modulations of the form described above. Furthermore, since the time-derivative energy $\ETan(t)\in L^1(0,\infty)$ is nonincreasing (see Proposition~\ref{lem:tangent_closed} and Lemma~\ref{lem:tangent_decay}), every global solution in Theorem~\ref{thm:main} also exhibits the refined asymptotic decay}
\[
\norm{f(t)}_{H^3}+\norm{Z(t)}_{H^3}=o(t^{-1/4}),\qquad
\norm{f_t(t)}_{\Hm}+\norm{Z_t(t)}_{H^1}=o(t^{-1/2})
\qquad\text{as }t\to\infty.
\]
The stability of left acoustic shocks follows from the reflection $x_1\mapsto-x_1$ and $u_1\mapsto-u_1$.
\end{remark}

Compared with the isentropic problem on $\R\times\mathbb T^2$~\cite{ww} and the constant-state problem on $\R^3$~\cite{gw}, proving Theorem~\ref{thm:main} for the full system \eqref{eq:nsf} in $\R^3$ involves several new difficulties.

\medskip\noindent\textbf{Difficulties:}
\begin{itemize}\setlength{\itemsep}{3pt}\setlength{\parskip}{0pt}
\item \textbf{Absence of a transverse Poincar\'e inequality on $\R^3$.}
{On $\R\times\mathbb T^2$, the transverse Poincar\'e inequality controls the non-zero transverse Fourier modes by their transverse derivatives, leaving the transverse average (the zero Fourier mode) as a one-dimensional problem. On $\R^3$, however, the continuous transverse spectrum extends to zero frequency, and integrating a conservation law $\partial_t\rho+\partial_1 m_1+\diver_\perp m_\perp=0$ in $x_1$ leaves the unintegrated transverse flux $\diver_\perp m_\perp$.}

\item \textbf{Opposite compression signs in $L^2(\R^3)$ and $\Hm(\R^3)$.}
In the standard $L^2(\R^3)$ energy estimate for $Z=(\phi,\psi,\zeta)^T$, the shock compression $\bu'<0$ generates an unfavorable positive term $\frac12\int|\bu'||\psi_1|^2\dd$ on the right-hand side, which cannot be bounded by $\norm{\nabla(\psi,\zeta)}_2^2$ on $\R^3$. {Passing to the unweighted $\Hm(\R^3)$ norm~\cite{gw} converts the normal derivative into a dissipative term, but simultaneously reverses the transverse gradient into an unfavorable positive term $+\frac12\int|\bu'||\nabla_\perp h|^2\dd$ (see \eqref{eq:model-scalar-inverse}).}

\item \textbf{Lack of density dissipation and coupling with temperature.}
In \eqref{eq:nsf}, temperature is coupled to density and velocity through the pressure law $p=\rho\theta$ and viscous heating $S(u):\nabla u$. Because the continuity equation lacks viscosity, a direct $\Hm(\R^3)$ estimate on the conservative variables fails to control $D^2h_\rho$, while the spatially varying shock profile couples the acoustic, entropy, and vorticity modes at first order.
\end{itemize}

\medskip\noindent\textbf{Strategies:}
\begin{itemize}\setlength{\itemsep}{3pt}\setlength{\parskip}{0pt}
\item \textbf{Entropy-weighted $\Hm(\R^3)$ multiplier and modified curl operator.}
To extract a favorable sign from $\bu'<0$ in both normal and transverse directions, we introduce the elliptic operator $\cA:\Hp(\R^3;\R^5)\to\Hm(\R^3;\R^5)$,
\[
\cA h=-\sum_{j=1}^3\partial_j(\omega M\partial_jh)+J^T\curl\bigl(\omega\br\bt\,\curl(Jh)\bigr),
\]
where $M=\eta_{UU}(\bar U)^{-1}$ is the inverse entropy Hessian, $\omega=(\br/\rho_*)^k$ is a density weight with a large exponent $k>0$, and $Jh=h_m+\bu h_E e_1$. {In the base Euler pairing with $h=\cA^{-1}g$, the symmetry of $A_jM$ reduces all second-order derivatives of $h$ to exact divergences, leaving a quadratic form in $\nabla h$ weighted by profile derivatives. In the rotational pairing, the alignment $\nabla\bu\times e_1=0$ ensures that $\Omega_h=\curl(Jh)$ involves only first-order derivatives of $h$; taking the curl eliminates the acoustic gradients, and the remaining normal transport term integrates by parts into a first-order quadratic form in \eqref{eq:idea_neg_curl_pairing}.}

\item \textbf{Cancellation on the acoustic kernel via shock compression.}
{By relative mass conservation, the weight derivative satisfies $\omega'/\omega=k\bu'/(s-\bu)<0$. At the reference state, factoring out $\omega\bu'$ from the weight-derivative terms in the entropy pairing yields a nonnegative quadratic form whose kernel consists of gradients $\partial_j h$ aligned with the left acoustic eigenvector $\ell$. Adding the rotational curl weight controls the two transverse acoustic components $(\partial_2h,\partial_3h)$ and reduces the joint null space to the 1-D normal acoustic line $\nabla h=(\alpha\ell,0,0)$. On this normal line, in the normalized units of \eqref{eq:idea_neg_tangent}, the left eigenvector identity $\ell^T(A_{1*}-cI)=0$ eliminates the metric-derivative term $\ell^T(A_{1*}-cI)M_*'\ell$, while differentiating the eigenvalue relation $A_1r=\lambda_+r$ with $r=M_*\ell$ and $\ell^Tr=2$ gives}
\begingroup\[
(\ell^TA_1'r)_*=2(\lambda_+')_*=(\gamma+1)\bu'.
\]
\endgroup
{Since the characteristic speed is compressive along the shock profile, $\lambda_+'=((\gamma+1)/2+O(\delta))\bu'<0$ by \eqref{eq:idea_neg_actual_speed},}
choosing $k>0$ sufficiently large and then $\delta>0$ small renders the full Euler pairing strictly negative definite:
\[
Q_E\le -c_E\int\omega|\bu'|\,|\nabla h|^2\dd.
\]

\item \textbf{Effective momentum and density dissipation estimate.}
We introduce the effective momentum $n=\mm+d\nabla\rr$ and evaluate the $\Hm(\R^3)$ energy on $g=(\rr,n,\EE)^T$. This creates an artificial mass diffusion $d\Delta\rr$ in the continuity equation so that the constant-state diffusive pairing controls $\norm{D^2h}_2^2$, while third-order density derivatives in the momentum equation are integrated by parts onto $D^2h$ at the expense of $\norm{\nabla\rr}_2^2$. We then bound $\norm{\nabla\phi}_2^2$ via the momentum--density cross functional $\int\psi\cdot\nabla\phi\dd$ and absorb the unfavorable compression term from the $L^2$ energy.

\item \textbf{Higher-order energy closure and time-derivative decay.}
{At top order ($|\alpha|=3$), integrating the pressure coupling $(\theta-\theta_*)\psi_\alpha\cdot\nabla\phi_\alpha$ by parts pairs $\phi_\alpha$ with $(\theta-\theta_*)\diver\psi_\alpha+\psi_\alpha\cdot\nabla\theta$; the first term is controlled by the fourth-order velocity dissipation and the second by $\norm{\nabla\theta}_\infty$, so no fourth-order density derivative is needed.} Finally, since the profile $\bar U(x_1)$ is stationary in the moving frame, the time-differentiated perturbation $(\dot f,\dot Z)=(f_t,Z_t)$ inherits the same linear structure, yielding ${\ETan}(t)\le C\varepsilon^2(1+t)^{-1}$ and the decay rates \eqref{eq:rate}--\eqref{eq:trate}.
\end{itemize}

\noindent\textbf{Notations:} Throughout this paper, we use the following notations.
\begin{itemize}\setlength{\itemsep}{2pt}\setlength{\parskip}{0pt}
\item All spatial norms are taken over $\R^3$. We denote by $\norm{\cdot}_p$ the $L^p(\R^3)$ norm ($1\le p\le\infty$) and by $H^k(\R^3)$ the usual Sobolev space of order $k\ge0$.
\item Using the Fourier transform $\widehat v(\xi)$, the homogeneous negative Sobolev space $\Hm(\R^3)$ is defined by
\[
\norm v_{\Hm}^2=\int_{\R^3}|\xi|^{-2}|\widehat v(\xi)|^2\,d\xi,
\]
and $\Hp(\R^3)$ denotes the completion of $C_c^\infty(\R^3)$ under $\norm{\nabla\cdot}_2$, with $\Hp(\R^3)\hookrightarrow L^6(\R^3)$.
\item We write $x_\perp=(x_2,x_3)$ and $\nabla_\perp=(\partial_2,\partial_3)$ for the transverse coordinates and gradient, and $D^2h$ for the full spatial Hessian.
\item Generic positive constants $c,C>0$ depend only on $(\mu,\lambda,\kappa,c_v)$ and the reference-state set, and are independent of $\delta$ and $\varepsilon$. We write $A\asymp B$ if $cB\le A\le CB$.
\end{itemize}

\leavevmode{The rest of the paper is organized as follows. Section~\ref{sec:prelim} presents the main analytical mechanism, shock profile properties, local well-posedness, and elliptic regularity of the weighted multiplier $\cA$. Section~\ref{sec:compression} establishes the weighted $\Hm(\R^3)$ energy inequality and the $H^3(\R^3)$ energy estimates in physical variables. Section~\ref{sec:global_decay} closes the coupled a priori estimates and proves global existence and the time-decay rates.}

\section[Preliminaries]{{Preliminaries}}\label{sec:prelim}

\leavevmode{We first examine a 3-D scalar viscous Burgers shock to illustrate why passing from $L^2(\R^3)$ to an unweighted $\Hm(\R^3)$ norm reverses the normal compression sign but fails in the transverse directions.} We then present the main algebraic cancellations in our weighted $\Hm(\R^3)$ multiplier and outline how the coupled energy estimates close for the full system \eqref{eq:nsf}.

Consider a planar Burgers shock $\bar u(x_1-st)$ in $\R^3$ satisfying
\[
\bar u''=(\bar u-s)\bar u',\qquad \bar u'<0.
\]
In the moving frame, the perturbation $v=u-\bar u$ solves
\begin{equation}\label{eq:model-burgers}
v_t+\partial_1\left\{(\bar u-s)v+\frac12v^2\right\}=\Delta v.
\end{equation}
Multiplying \eqref{eq:model-burgers} by $v$ and integrating by parts over $\R^3$ gives the $L^2$ energy identity
\begin{equation}\label{eq:model-scalar-physical}
\frac12\frac{d}{dt}\norm v_2^2+\norm{\nabla v}_2^2
=-\frac12\int\bar u'v^2\dd
=\frac12\int|\bar u'|\,v^2\dd.
\end{equation}
\leavevmode{For a nontrivial shock profile on $\R^3$, the compression integral on the right-hand side of \eqref{eq:model-scalar-physical} cannot be controlled by $\norm{\nabla v}_2^2$: testing with a rescaled cutoff $v_R(x)=\chi(x_1/R)\eta(x_\perp/R)$ satisfying $\chi(0)=\eta(0)=1$, the compression integral $\int|\bar u'|v_R^2\dd$ grows like $R^2$ as $R\to\infty$ since $\int_{\R}|\bar u'|\,dx_1>0$, whereas $\norm{\nabla v_R}_2^2$ grows only like $R$.}

\leavevmode{For smooth perturbations in $H^1(\R^3)\cap\Hm(\R^3)$, testing \eqref{eq:model-burgers} instead with the unweighted potential $h=(-\Delta)^{-1}v\in\Hp(\R^3)$ gives $\ip{\Delta v}h=-\norm v_2^2$, while integrating the drift term by parts yields}
\[
\int(\bar u-s)(-\Delta h)\partial_1h\dd
=\frac12\int\bar u'\bigl(|\partial_1h|^2-|\nabla_\perp h|^2\bigr)\dd.
\]
Thus the unweighted $\Hm(\R^3)$ energy satisfies
\begin{equation}\label{eq:model-scalar-inverse}
\frac12\frac{d}{dt}\norm{\nabla h}_2^2+\norm v_2^2
+\frac12\int|\bar u'|\,|\partial_1h|^2\dd
=\frac12\int|\bar u'|\,|\nabla_\perp h|^2\dd
+\frac12\int v^2\partial_1h\dd,
\end{equation}
where the cubic term is bounded by $C\norm v_{H^1}\norm v_2^2$ via $\norm{\nabla h}_6\le C\norm v_2$. We see from \eqref{eq:model-scalar-inverse} that passing to $\Hm(\R^3)$ makes the normal derivative $\partial_1h$ dissipative on the left-hand side, but flips the transverse gradient $\nabla_\perp h$ to the positive sign on the right-hand side.

For the full system \eqref{eq:nsf}, let $f=(\rr,\mm,\EE)^T$ and $Z=(\phi,\psi,\zeta)^T$ be defined by \eqref{eq:variables}, and set $N=\norm Z_{H^3}$. We introduce the effective momentum $n$ and modified conservative state $g$:
\begin{gather}
n=\mm+d\nabla\rr,\qquad g=(\rr,n,\EE)^T,\nonumber\\
d=\min\left\{\frac{\nu}{2\rho_*},\,\frac{\kappa}{\rho_*}\right\}>0.
\label{eq:idea_neg_effective}
\end{gather}
Let $\eta(U)=-\rho(c_v\log\theta-\log\rho)$ be the physical entropy and $M=\eta_{UU}(\bar U)^{-1}$ be its inverse Hessian. For a potential $h=(h_\rho,h_m,h_E)^T$ with $h_m=(h_{m,1},h_{m,2},h_{m,3})^T\in\R^3$, we define $J$, $J^T$, $\Omega_h$, and the elliptic operator $\cA$ by
\begin{gather}
Jh=h_m+\bu h_Ee_1,\qquad
J^Tz=(0,z,\bu z_1)^T,\qquad \Omega_h=\curl(Jh),\nonumber\\
\omega=\left(\frac{\br}{\rho_*}\right)^k,\qquad
g=\cA h:=-\sum_{j=1}^3\partial_j(\omega M\partial_jh)
+J^T\curl(\omega\br\bt\Omega_h),
\label{eq:idea_neg_multiplier}
\end{gather}
where $k>0$ is a large constant to be chosen independently of $\delta$. By Lemma~\ref{lem:inverse}, each $g\in\Hm(\R^3)$ determines a unique $h=\cA^{-1}g\in\Hp(\R^3)$, and the negative-order energy is
\begin{equation}\label{eq:idea_neg_energy}
\cE_-(g)=\frac12\ip gh
=\frac12\int\omega\left\{\sum_{j=1}^3(\partial_jh)^TM\partial_jh
+\br\bt|\Omega_h|^2\right\}\dd
\asymp\norm g_{\Hm}^2.
\end{equation}
Our key step is to show that the linearized Euler contribution $Q_E$ in $\cE_-'(g)$ satisfies
\begin{equation}\label{eq:idea_neg_target}
Q_E\le-cD_s,\qquad
D_s:=\int\omega|\bu'|\,|\nabla h|^2\dd.
\end{equation}

\medskip\noindent\textbf{1. Divergence reduction of the entropy pairing.}
Let $A_j=F_{j,U}(\bar U)$ and $\widetilde A_j=A_j-s\delta_{j1}I$. Taking the inner product of the effective perturbation equation with $h=\cA^{-1}g$ gives
\begin{equation}\label{eq:idea_neg_test}
\begin{gathered}
g_t+\sum_{j=1}^3\partial_j(\widetilde A_jg)=\mathcal R,
\qquad Y=\sum_{j=1}^3\widetilde A_j^T\partial_jh,\\
\cE_-'=\ip{g_t}h=Q_E+\ip{\mathcal R}h,
\qquad Q_E=\ip{\cA h}Y,
\end{gathered}
\end{equation}
where $\mathcal R=\mathcal B g+\cN_g$ collects the diffusive and nonlinear terms in \eqref{eq:rhs} and \eqref{eq:effective_nonlinearity}. By entropy symmetrization, the matrices $B_j:=\widetilde A_jM=M\widetilde A_j^T$ are symmetric. Expanding the entropy-metric part $Q_M$ of $Q_E=\ip{\cA h}Y$ gives
\begin{align*}
Q_M
&=\int\omega\sum_{i,j=1}^3(\partial_ih)^TM
\partial_i(\widetilde A_j^T\partial_jh)\dd\\
&=\int\omega\sum_{j=1}^3(\partial_1h)^TM(A_j')^T\partial_jh\dd
+\int\omega\sum_{i,j=1}^3
(\partial_ih)^TB_j\partial_{ji}h\dd.
\end{align*}
Using the symmetry of $B_j$ to write
\[
(\partial_ih)^TB_j\partial_{ji}h
=\frac12\partial_j\bigl\{(\partial_ih)^TB_j\partial_ih\bigr\}
-\frac12\delta_{j1}(\partial_ih)^TB_1'\partial_ih
\]
and integrating by parts over $\R^3$, we obtain
\begin{equation}\label{eq:idea_neg_entropy_cancel}
Q_M
=-\frac12\int\sum_{i=1}^3 (\partial_ih)^T
(\omega\widetilde A_1M)'\partial_ih\dd
+\int\omega\sum_{j=1}^3(\partial_1h)^TM(A_j')^T\partial_jh\dd.
\end{equation}
Hence all second derivatives $D^2h$ cancel out of $Q_M$.

\medskip\noindent\textbf{2. Elimination of second-order gradients in the curl pairing.}
In the curl part, since $\nabla\bu=\bu'e_1$ is parallel to $e_1$, we have $\nabla\bu\times e_1=0$, so that
\begin{equation}\label{eq:idea_neg_no_zeroth}
\curl(\bu h_Ee_1)
=\bu\nabla h_E\times e_1.
\end{equation}
Thus $\Omega_h=\curl(Jh)$ involves only first derivatives of $h$. With $i_0=(c_v+1)\bt+\bu^2/2$, direct calculation gives
\begin{align}
JY={}&\nabla h_\rho+\bu\nabla h_{m,1}+i_0\nabla h_E
+(\bu-s)(\partial_1h_m+\bu e_1\partial_1h_E),\nonumber\\
\curl(JY)={}&(\bu-s)\partial_1\Omega_h
+\bu'e_1\times(\nabla h_{m,1}+\partial_1h_m+\bu\nabla h_E)\nonumber\\
&+(i_0'-s\bu')e_1\times\nabla h_E.
\label{eq:idea_neg_curl_cancel}
\end{align}
\leavevmode{Here the curl operator annihilates the acoustic gradients $\nabla h_\rho$, $\nabla h_{m,1}$, and $\nabla h_E$, while pairing the remaining second-order transport term $(\bu-s)\partial_1\Omega_h$ with $\Omega_h$ and integrating by parts gives}
\begin{align}
Q_\Omega&=\int\omega\br\bt\,\Omega_h\cdot\curl(JY)\dd,
\qquad Q_E=Q_M+Q_\Omega,\nonumber\\
\int\omega\br\bt(\bu-s)\Omega_h\cdot\partial_1\Omega_h\dd
&=-\frac12\int[(\bu-s)\omega\br\bt]'|\Omega_h|^2\dd.
\label{eq:idea_neg_curl_pairing}
\end{align}
Therefore, $Q_E=Q_M+Q_\Omega$ reduces to a quadratic form in $\nabla h$ weighted by profile derivatives.

\medskip\noindent\textbf{3. Coercivity on the acoustic kernel via shock compression.}
From the relative mass conservation $[\br(\bu-s)]'=0$, the derivative of $\omega$ shares the negative sign of $\bu'<0$:
\begin{equation}\label{eq:idea_neg_mass_weight}
[\br(\bu-s)]'=0
\quad\Longrightarrow\quad
\frac{\br'}{\br}=\frac{\bu'}{s-\bu},
\qquad \frac{\omega'}{\omega}=k\frac{\bu'}{s-\bu}<0.
\end{equation}
\leavevmode{By the natural scaling transformation $U=\rho_*T\tilde U$, $h=T^{-1}\tilde h$, and $s=q\tilde s$ with $q=\sqrt{\theta_*}$ and $T=\diag(1,q,q,q,q^2)$, the reference state $(\rho_*,0,\theta_*)$ is normalized to $(1,0,1)$, while the matrices in \eqref{eq:Eulerflux} and \eqref{eq:M} and the operator $J$ transform as}
\begingroup\[
M=\rho_*T\tilde M T,\qquad
A_j=qT\tilde A_jT^{-1},\qquad
Jh=q^{-1}\tilde J\tilde h.
\]
\endgroup
\leavevmode{Under this rescaling, the weight $\omega=(\br/\rho_*)^k$ is invariant, $\curl(JY)=\curl(\tilde J\tilde Y)$, and both $Q_M$ and $Q_\Omega$ in \eqref{eq:idea_neg_entropy_cancel} and \eqref{eq:idea_neg_curl_pairing} scale by the common positive factor $\rho_*q$. Since $\rho_*$, $\theta_*$, and the singular values of $T$ are uniformly bounded above and away from zero on $K$, it suffices to verify \eqref{eq:idea_neg_target} in normalized units where $\rho_*=\theta_*=1$ and $c=\sqrt\gamma$.}

\leavevmode{For $\rho_*=\theta_*=1$ and $c=\sqrt\gamma$, Lemma~\ref{lem:profile} yields}
\begin{gather}
|\br-1|+|\bt-1|+|\bu|+|s-c|\le C\delta,\nonumber\\
\frac{\br'}{\bu'}=\frac1c+O(\delta),\qquad
\frac{\bt'}{\bu'}=\frac{\gamma-1}{c}+O(\delta),\qquad
\frac{\omega'}{\omega\bu'}=\frac kc+O_k(\delta).
\label{eq:idea_neg_profile_comparison}
\end{gather}
\leavevmode{Consequently, along the shock profile, the spatial derivative of any smooth matrix $B(\rho,u_1,\theta)$ satisfies}
\begingroup\[
B'(\bar U)=\br'B_\rho(\bar U)+\bu'B_{u_1}(\bar U)+\bt'B_\theta(\bar U)
=B_*'+O(\delta|\bu'|),
\]
\endgroup
\leavevmode{where $B_*'=\bu'(x_1)\bigl(c^{-1}\partial_\rho+\partial_{u_1}+(\gamma-1)c^{-1}\partial_\theta\bigr)B|_{(1,0,1)}$ denotes the leading-order directional derivative at the reference state $(1,0,1)$ along the acoustic tangent, corresponding to the substitutions}
\begin{equation}\label{eq:idea_neg_tangent}
\br=\bt=1,\quad\bu=0,\quad s=c,\qquad
\br'=\frac{\bu'}{c},\quad\bt'=\frac{\gamma-1}{c}\bu',
\quad\frac{\omega'}{\omega}=\frac{k}{c}\bu'.
\end{equation}
Denoting evaluation at the reference state by a subscript $*$, the weight-derivative term in \eqref{eq:idea_neg_entropy_cancel} is
\begin{equation}\label{eq:idea_neg_weight}
-\frac{\omega'}{2}(\partial_jh)^T(A_{1*}-cI)M_*\partial_jh
=\frac k2\omega\bu'(\partial_jh)^T(M_*-c^{-1}A_{1*}M_*)\partial_jh.
\end{equation}
For any $X=(a,b_1,b_2,b_3,e)^T\in\R^5$, using $c^2=(c_v+1)/c_v$, we complete the square to obtain
\begin{align}
\Phi(X)&:=X^T(M_*-c^{-1}A_{1*}M_*)X\nonumber\\
&=a^2+2c_v a e+c_v(c_v+1)e^2+|b|^2
-\frac2c[a+(c_v+1)e]b_1\nonumber\\
&=\left(\frac{a+(c_v+1)e}{c}-b_1\right)^2
+\frac{a^2}{c_v+1}+b_2^2+b_3^2.
\label{eq:idea_neg_phi}
\end{align}
Hence $\Phi(X)\ge0$, and its null space is spanned by the left acoustic eigenvector $\ell$:
\[
\sum_{j=1}^3\Phi(\partial_jh)=0
\quad\Longleftrightarrow\quad
\partial_jh=\alpha_j\ell\quad(j=1,2,3),\qquad
\ell=\bigl(0,\,1,\,0,\,0,\,(c_vc)^{-1}\bigr)^T.
\]
At the reference state where $\bu=0$, we have $\Omega_h=\curl h_m$, so the curl weight term in \eqref{eq:idea_neg_curl_pairing} adds $|\curl h_m|^2$. The combined weight-derivative terms equal $\frac k2\omega\bu'F(\nabla h)$, where
\begin{equation}\label{eq:idea_neg_square}
F(\nabla h)=\sum_{j=1}^3\Phi(\partial_jh)+|\curl h_m|^2\ge0.
\end{equation}
When $\partial_jh=\alpha_j\ell$ (so $\partial_jh_m=\alpha_je_1$), we have
\begin{gather}
\curl h_m=\alpha\times e_1=(0,\alpha_3,-\alpha_2)^T,\nonumber\\
\ker F=\{(\alpha\ell,0,0):\alpha\in\R\}.
\label{eq:idea_neg_kernel}
\end{gather}
Thus the curl term supplies $\alpha_2^2+\alpha_3^2$ and reduces $\ker F$ to the 1-D normal acoustic line.

Let $P(\nabla h)$ denote the quadratic form at the reference state coming from the terms in $Q_M+Q_\Omega$ where the derivative hits the profile matrices rather than $\omega$. Then $P$ is independent of $k$, and the leading part of $Q_E$ is
\begin{equation}\label{eq:idea_neg_tangent_pairing}
Q_{E,*}=\int\omega\bu'
\left\{P(\nabla h)+\frac k2F(\nabla h)\right\}\dd.
\end{equation}
Set
\begin{equation}\label{eq:idea_neg_eigenvectors}
r=M_*\ell=\bigl(c^{-1},\,1,\,0,\,0,\,c_vc\bigr)^T,\qquad
A_{1*}r=cr,\quad\ell^TA_{1*}=c\ell^T,\quad\ell^Tr=2.
\end{equation}
For $\nabla h=(\alpha\ell,0,0)\in\ker F$, all curl terms in $P$ vanish, and the two terms in \eqref{eq:idea_neg_entropy_cancel} give
\begin{align}
\omega\bu'P(\alpha\ell,0,0)
&=\frac{\omega\alpha^2}{2}
\left\{\ell^TA_1'M_*\ell
-\ell^T(A_{1*}-cI)M'\ell\right\}\nonumber\\
&=\frac{\omega\alpha^2}{2}\ell^TA_1'r,
\label{eq:idea_neg_metric_cancel}
\end{align}
since $\ell^T(A_{1*}-cI)=0$. {Extending $r$ smoothly in a neighborhood of the reference state as a right eigenvector of $A_1$ for the simple eigenvalue $\lambda_+$, differentiating $A_1r=\lambda_+r$ along the acoustic tangent \eqref{eq:idea_neg_tangent}, and multiplying on the left by $\ell^T$ gives}
\begin{align}
A_1'r+A_1r'&=\lambda_+'r+\lambda_+r',\nonumber\\
\ell^TA_1'r+c\ell^Tr'&=2\lambda_+'+c\ell^Tr',\nonumber\\
(\ell^TA_1'r)_*
&=2(\lambda_+')_*
=(\gamma+1)\bu',
\label{eq:idea_neg_eigenvector_cancel}
\end{align}
where $c\ell^Tr'$ cancels on both sides. Since $\lambda_+=\bu+\sqrt{\gamma\bt}$ satisfies
\begin{equation}\label{eq:idea_neg_actual_speed}
\lambda_+'=\bu'+\frac{\sqrt\gamma}{2\sqrt{\bt}}\bt'
=\left(\frac{\gamma+1}{2}+O(\delta)\right)\bu'<0,
\end{equation}
it holds on $\ker F$ that
\begin{equation}\label{eq:idea_neg_normal_sign}
P(\alpha\ell,0,0)=\frac{\gamma+1}{2}\alpha^2>0
\qquad(\alpha\ne0).
\end{equation}
Splitting $\nabla h=v+w$ with $v\in\ker F$ and $w\in(\ker F)^\perp$, we have for constants $a_0,b_0,C>0$ independent of $k$ and $\delta$ that
\begin{gather}
P(v)\ge a_0|v|^2,\quad F(\nabla h)=F(w)\ge b_0|w|^2,\quad
|P(v,w)|\le C|v||w|,\quad |P(w)|\le C|w|^2,\nonumber\\
P(\nabla h)+\frac k2F(\nabla h)
\ge\frac{a_0}{2}|v|^2+
\left(\frac{kb_0}{2}-C-\frac{2C^2}{a_0}\right)|w|^2
\ge c_0|\nabla h|^2,
\label{eq:idea_neg_coercivity}
\end{gather}
{for $k>0$ chosen sufficiently large. Here $a_0,b_0,c_0>0$ depend only on $c_v$, because $F$ and $P$ are strictly positive on the unit spheres of $(\ker F)^\perp$ and $\ker F$, respectively. With $k$ thus fixed, \eqref{eq:idea_neg_profile_comparison} implies that the actual profile coefficients in $Q_E/(\omega\bu')$ differ from their reference-state counterparts in $P(\nabla h)+\frac k2F(\nabla h)$ by at most $O_k(\delta)$, uniformly on the compact set $K$.} Since $\bu'<0$, choosing $\delta>0$ sufficiently small after fixing $k$ gives
\begin{equation}\label{eq:idea_neg_profile_sign}
Q_E\le-(c_0-C_k\delta)D_s\le-\frac{c_0}{2}D_s.
\end{equation}
\leavevmode{Returning to physical units, in the diffusive pairing $\ip{\mathcal Bg}h$, the longitudinal and rotational parts of $h_m$ are orthogonal via}
\begin{equation}\label{eq:idea_neg_divcurl}
\int\nabla\diver h_m\cdot\curl\curl h_m\dd
=-\int\diver h_m\,\diver(\curl\curl h_m)\dd=0,
\end{equation}
and the choice of $d$ in \eqref{eq:idea_neg_effective} makes the $(\rho,E)$ block coercive in $D_0=\norm{D^2h}_2^2$. Bounding the variable-coefficient and nonlinear terms (Proposition~\ref{prop:negative_energy}) yields
\begin{equation}\label{eq:idea_neg_full_estimate}
\cE_-'+c(D_0+D_s)
\le C\norm{\nabla\rr}_2^2
+CN\left\{D_0+\norm Z_2^2+
\norm{\nabla(\psi,\zeta)}_2^2\right\}.
\end{equation}

We now summarize how \eqref{eq:idea_neg_full_estimate} is combined with the energy estimates for $Z=(\phi,\psi,\zeta)^T$.

\medskip\noindent\textbf{1. Compression error in the $L^2$ energy.}
With $a=\rho^{-1}$, the system for $Z=(\phi,\psi,\zeta)^T$ reads
\[
\left\{
\begin{aligned}
&D_t\phi+\diver\psi=F_\phi,\\
&D_t\psi+\theta\nabla\phi+\nabla\zeta-a\cL\psi=F_\psi,\\
&c_vD_t\zeta+\theta\diver\psi-\kappa a\Delta\zeta=F_\zeta,
\end{aligned}
\right.
\]
where $F_\phi,F_\psi,F_\zeta$ are given in \eqref{eq:physical_sources} and $F_\psi$ contains $-\bu'\psi_1e_1$. Multiplying the three equations by $\theta_*\phi$, $\psi$, and $\zeta/\theta_*$ defines
\[
S_0(Z)=\frac12\int\left(\theta_*\phi^2+|\psi|^2+
\frac{c_v}{\theta_*}\zeta^2\right)\dd\asymp\norm Z_2^2.
\]
The constant-coefficient pressure and temperature couplings integrate to zero:
\[
\int\theta_*(\phi\diver\psi+\psi\cdot\nabla\phi)\dd=0,
\qquad
\int(\psi\cdot\nabla\zeta+\zeta\diver\psi)\dd=0,
\]
whereas the normal velocity drift from transport integration by parts and $-\bu'\psi_1e_1$ yields the positive term
\begin{equation}\label{eq:intro_comp}
\begin{split}
R_0
&=\frac12\int\bu'|\psi_1|^2\dd-\int\bu'|\psi_1|^2\dd\\
&=\frac12\int|\bu'|\,|\psi_1|^2\dd
\le C\delta^2\norm Z_2^2,
\end{split}
\end{equation}
leading to the $L^2$ inequality (Lemma~\ref{lem:physical_l2})
\begin{equation}\label{eq:intro_l2}
S_0(Z)'+c\norm{\nabla(\psi,\zeta)}_2^2
\le C(\delta+N)\norm Z_2^2.
\end{equation}

\medskip\noindent\textbf{2. Elliptic recovery and density dissipation.}
Set $D_1(Z)=\norm{\nabla\phi}_2^2+\norm{\nabla(\psi,\zeta)}_{H^1}^2$. From the relations
\begin{gather*}
\rr=\br(e^\phi-1),\qquad
\nabla\rr=\rho\nabla\phi+\br'(e^\phi-1)e_1,\qquad
\rho\psi=n-\bar u\rr-d\nabla\rr,\\
c_v\rho\zeta=\EE-
\rr\left(c_v\bt+\frac12|\bar u|^2\right)
-\rho\bar u\cdot\psi-\frac12\rho|\psi|^2,
\end{gather*}
we have
\begin{equation}\label{eq:intro_recover}
\norm Z_2^2\le C\bigl(\norm g_2^2+\norm{\nabla\phi}_2^2\bigr),\qquad
R_0\le C\delta^2\bigl(\norm g_2^2+\norm{\nabla\phi}_2^2\bigr),
\end{equation}
and by Lemma~\ref{lem:elliptic},
\begin{equation}\label{eq:idea_neg_recovery_bridge}
\norm g_2^2\le C(D_0+\delta^2D_s),\qquad
\norm{\nabla\rr}_2^2\le C\norm{\nabla\phi}_2^2+C\delta^4\norm g_2^2.
\end{equation}
Substituting \eqref{eq:intro_recover}--\eqref{eq:idea_neg_recovery_bridge} into \eqref{eq:idea_neg_full_estimate} gives
\begin{equation}\label{eq:intro_inverse}
\cE_-'+c_g\norm g_2^2\le C_*D_1(Z).
\end{equation}
To control $\norm{\nabla\phi}_2^2$ on the right-hand side of \eqref{eq:intro_inverse}, we use the cross functional $I_1(Z)=\int\psi\cdot\nabla\phi\dd$, which yields (Lemma~\ref{lem:physical_h1})
\[
I_1(Z)'+c_I\norm{\nabla\phi}_2^2\le C_I\norm{\nabla(\psi,\zeta)}_{H^1}^2+C(\delta+N)\norm Z_{H^1}^2.
\]
Setting $P_1(Z)=S_1(Z)+\sigma_1 I_1(Z)\asymp\norm Z_{H^1}^2$ for small $\sigma_1>0$, we obtain
\begin{equation}\label{eq:intro_physical}
P_1(Z)'+c_pD_1(Z)\le C_p(\delta+N)\norm Z_2^2
\le C_p'(\delta+N)\bigl\{\norm g_2^2+D_1(Z)\bigr\}.
\end{equation}
Choosing $L_1\ge1$ with $L_1c_p-C_*\ge1$ and adding $L_1$ times \eqref{eq:intro_physical} to \eqref{eq:intro_inverse} gives, for small $\delta+N$,
\begin{equation}\label{eq:intro_closed}
(\cE_-+L_1P_1)'+\frac{c_g}{2}\norm g_2^2+\frac12D_1(Z)\le0.
\end{equation}

\medskip\noindent\textbf{3. Higher-order closure and time-derivative decay.}
For $m=2,3$, we combine $S_m(Z)$ with $I_m(Z)=\sum_{|\alpha|\le m-1}\int\partial^\alpha\psi\cdot\nabla\partial^\alpha\phi\dd$ and integrate the pressure error by parts:
\[
-\int(\theta-\theta_*)\psi_\alpha\cdot\nabla\phi_\alpha\dd
=\int(\theta-\theta_*)\phi_\alpha\diver\psi_\alpha\dd
+\int\phi_\alpha\psi_\alpha\cdot\nabla\theta\dd.
\]
\leavevmode{For $|\alpha|=3$, this integral is bounded in absolute value by}
\begingroup\[
\norm{\phi_\alpha}_2\bigl(\norm{\theta-\theta_*}_\infty\norm{\diver\psi_\alpha}_2
+\norm{\nabla\theta}_\infty\norm{\psi_\alpha}_2\bigr),
\]
\endgroup
{where the only fourth-order derivative falls on $\psi$ and is absorbed by the viscous dissipation.} {Finally, for the time-derivative energy $\ETan=\cE_-(g_t)+L_1P_1(Z_t)$, Proposition~\ref{lem:tangent_closed} and Lemma~\ref{lem:tangent_decay} show that $\ETan'\le0$ and}
\[
{\ETan}(0)+\int_0^\infty{\ETan}(t)\,dt\le C\varepsilon^2
\quad\Longrightarrow\quad
{\ETan}(t)\le C\varepsilon^2(1+t)^{-1},
\]
which combined with $|\cE_1'|\le 2\sqrt{\cE_1{\ETan}}$ and third-order parabolic damping yields \eqref{eq:rate}--\eqref{eq:trate}.

Let $z=y_1-st$ denote the traveling coordinate in the original frame, with a prime denoting $d/dz$. \leavevmode{Under the endstate assumptions of Theorem~\ref{thm:main}, the planar shock profile $\bar U(z)$ has constant negative relative mass flux $j=\br(\bu-s)<0$.} Integrating the one-dimensional profile equations from $-\infty$ to $z$ gives the first-order ODE system
\begin{equation}\label{eq:profile}
\left\{
\begin{aligned}
&\br=\frac{j}{\bu-s},\\
&\nu\bu'=j(\bu-u_{1,-})+\br\bt-\rho_-\theta_-,\\
&\kappa\bt'=jc_v(\bt-\theta_-)+\rho_-\theta_-(\bu-u_{1,-})
-\frac j2(\bu-u_{1,-})^2.
\end{aligned}
\right.
\end{equation}

\leavevmode{The standard small-amplitude profile estimates are summarized below (cf.~\cite[Lemma~2.3]{kvw}).}

\begin{lemma}\label{lem:profile}
For sufficiently small shock strength $\delta=\rho_--\rho_+>0$, the system \eqref{eq:profile} admits a smooth profile $\bar U(z)$ connecting $U_-$ to $U_+$ that is unique up to spatial translation within the small-amplitude class. Under the normalization $\bu(0)=0$ and \eqref{eq:reference}, the profile satisfies $\bu'(z)<0$ for all $z\in\R$ and
\begin{gather}
|\br-\rho_*|+|\bt-\theta_*|+|\bu|+\bigl|s-\sqrt{\gamma\theta_*}\bigr|\le C\delta,
\quad |\bu'(z)|\le C\delta^2e^{-c\delta|z|},\label{eq:prof1}\\
|\br'|+|\bt'|\le C|\bu'|,\qquad
|\br''|+|\bu''|+|\bt''|\le C\delta|\bu'|,\label{eq:prof2}\\
\norm{\partial_z^\ell(\br,\bu,\bt)}_\infty\le C_\ell\delta^2,
\qquad 1\le\ell\le5,\label{eq:prof3}\\
\frac{\br'}{\br}=\frac{\bu'}{s-\bu},\qquad
\left|\frac{\bt'}{\bt}-(\gamma-1)\frac{\br'}{\br}\right|
\le C\delta|\bu'|.\label{eq:prof4}
\end{gather}
\end{lemma}

\leavevmode{Since $\br$ and $\bt$ are uniformly bounded above and away from zero with bounded derivatives of all orders, the change of variables between $f=(\rr,\mm,\EE)^T$ and $Z=(\phi,\psi,\zeta)^T$ in \eqref{eq:variables} defines a smooth diffeomorphism on any $H^3(\R^3)$ neighborhood of the origin where $\rho,\theta>0$, and standard Sobolev composition estimates yield}
\begingroup\begin{equation}\label{eq:positiveequiv}
\norm f_{H^j}\asymp\norm Z_{H^j},\qquad 0\le j\le3.
\end{equation}
\endgroup
\leavevmode{By the standard local well-posedness theory for symmetrizable hyperbolic--parabolic systems (cf.~\cite{serre}), for any initial perturbation $Z_0\in H^3(\R^3)$ with $\inf_{\R^3}\rho_0>0$ and $\inf_{\R^3}\theta_0>0$, there exists a local lifespan $T>0$, depending only on $\norm{Z_0}_{H^3}$ and the lower bounds of $(\rho_0,\theta_0)$, such that the Cauchy problem admits a unique strong solution $Z\in C([0,T];H^3(\R^3))\cap C^1([0,T];H^1(\R^3))$ with $\nabla(\psi,\zeta)\in L^2(0,T;H^3(\R^3))$. The solution depends continuously on the initial data in $H^3(\R^3)$, preserves higher Sobolev regularity, and extends in time as long as $\norm{Z(t)}_{H^3}$ remains bounded and $(\rho,\theta)$ stay uniformly bounded away from zero.}

\leavevmode{To establish global existence, we work on an arbitrary existence interval $[0,T]$ under the a priori assumption}
\begingroup\begin{equation}\label{eq:bootstrap}
N(t):=\norm{Z(t)}_{H^3}\le\varepsilon_1,\qquad 0<\delta\le\delta_1,
\end{equation}
\endgroup
{with small constants $\delta_1,\varepsilon_1>0$ to be determined in \eqref{eq:threshold}.}

In the moving shock frame, subtracting the stationary profile equations from \eqref{eq:nsf} gives the conservative perturbation system
\begin{equation}\label{eq:conservative_perturbation}
f_t-s\partial_1f+\sum_{j=1}^3\partial_j
\{F_j(\bar U+f)-F_j(\bar U)\}
=\begin{pmatrix}
0\\
\cL\psi\\
\kappa\Delta\zeta+\diver\{S(\bar u)\psi+S(\psi)\bar u+S(\psi)\psi\}
\end{pmatrix}.
\end{equation}
To create mass diffusion at negative Sobolev order, we fix $d>0$ and define the effective momentum $n$ and modified conservative perturbation $g$ by
\begin{equation}\label{eq:d}
d=\min\left\{\frac{\nu}{2\rho_*},\,\frac{\kappa}{\rho_*}\right\}>0,
\qquad n=\mm+d\nabla\rr,\qquad g=(\rr,n,\EE)^T.
\end{equation}
Under this substitution, the continuity equation becomes
\[
\rr_t-s\partial_1\rr+\diver n=d\Delta\rr,
\]
and the higher-order density derivatives introduced into the momentum and energy equations will be integrated by parts onto the elliptic potential.

Writing \eqref{eq:conservative_perturbation} in divergence form as $f_t=-\diver\mathcal F$, where $\mathcal F=(\mathcal F_1,\mathcal F_2,\mathcal F_3)$ denotes the total perturbation flux, we have on any local interval $[0,T]$ satisfying \eqref{eq:bootstrap} that
\[
\norm{f_t}_{\Hm}\le\norm{\mathcal F}_2\le C\bigl(\norm f_2+\norm{\nabla\psi}_2+\norm{\nabla\zeta}_2\bigr).
\]
Since the continuity equation also implies $\rr_t\in C([0,T];L^2(\R^3))$, the identity $\norm{\nabla\rr_t}_{\Hm}=\norm{\rr_t}_2$ yields
\[
g_t=f_t+(0,d\nabla\rr_t,0)^T\in L^2(0,T;\Hm(\R^3)),
\]
so that
\[
f,g\in C([0,T];\Hm(\R^3))\cap H^1(0,T;\Hm(\R^3))
\]
on every local existence interval whenever $f_0\in H^3(\R^3)\cap\Hm(\R^3)$.

Let $A_j=F_{j,U}(\bar U)$ denote the Euler flux Jacobians along the profile $\bar U$, and consider the physical entropy
\[
\eta(U)=-\rho(c_v\log\theta-\log\rho).
\]
\leavevmode{The entropy variable $\Lambda=\nabla_U\eta(U)$ has components}
\begingroup\[
\Lambda_\rho=-(c_v\log\theta-\log\rho)+c_v+1-\frac{|u|^2}{2\theta},
\qquad \Lambda_m=\frac{u}{\theta},
\qquad \Lambda_E=-\frac1\theta.
\]
\endgroup
Writing $e=c_v\theta+|u|^2/2$, its inverse differentials are given by
\[
d\theta=\theta^2d\Lambda_E,\qquad
du=\theta\,d\Lambda_m+\theta u\,d\Lambda_E,\qquad
\frac{d\rho}{\rho}=d\Lambda_\rho+u\cdot d\Lambda_m+e\,d\Lambda_E.
\]
Evaluating the Jacobian $\partial U/\partial\Lambda$ at $\bar U$ gives the inverse entropy Hessian $M=\eta_{UU}(\bar U)^{-1}$:
\begin{equation}\label{eq:M}
M=\br\begin{pmatrix}
1&\bar u^T&e_0\\
\bar u&\bt I+\bar u\bar u^T&i_0\bar u\\
e_0&i_0\bar u^T&e_0^2+c_v\bt^2+\bt\bu^2
\end{pmatrix},\qquad e_0=c_v\bt+\frac12\bu^2,\qquad i_0=e_0+\bt.
\end{equation}
Compatibility of the entropy flux $\eta(U)u_j$ with the Euler flux implies that $\eta_{UU}A_j$ is symmetric, and hence $A_jM=MA_j^T$. For any potential $h=(h_\rho,h_m^T,h_E)^T\in\R^5$ with $h_m=(h_{m,1},h_{m,2},h_{m,3})^T$, the quadratic form of $M$ decomposes into a sum of squares:
\[
h^TMh=\br\left\{(h_\rho+\bu h_{m,1}+e_0h_E)^2
+\bt|h_m+\bu h_Ee_1|^2+c_v\bt^2h_E^2\right\}.
\]

\leavevmode{With the spatial density weight $\omega(x_1)=(\br(x_1)/\rho_*)^k$ ($k>0$ fixed in Proposition~\ref{prop:negative_energy}), we define the linear map $J(x_1):\R^5\to\R^3$, its transpose $J^T(x_1):\R^3\to\R^5$, and the modified curl $\Omega_h$ by}
\[
Jh=h_m+\bu h_Ee_1,\qquad
J^Tz=(0,z,\bu z_1)^T\quad(z\in\R^3),\qquad
\Omega_h=\curl(Jh).
\]
For $h,\varphi\in\Hp(\R^3;\R^5)$, we define the stationary symmetric bilinear form
\begin{equation}\label{eq:form}
a(h,\varphi)=\int\omega\sum_{j=1}^3(\partial_jh)^TM\partial_j\varphi\dd
+\int\omega\br\bt\,\Omega_h\cdot\Omega_\varphi\dd.
\end{equation}

\begin{lemma}\label{lem:inverse}
The bilinear form \eqref{eq:form} defines a topological isomorphism $\cA:\Hp(\R^3;\R^5)\to\Hm(\R^3;\R^5)$ via $\ip{\cA h}\varphi=a(h,\varphi)$. For every $g\in\Hm(\R^3;\R^5)$, there exists a unique potential $h=\cA^{-1}g\in\Hp(\R^3;\R^5)$, uniquely represented in $L^6(\R^3;\R^5)$, such that $a(h,\varphi)=\ip g\varphi$ for all $\varphi\in\Hp(\R^3;\R^5)$, where in the sense of distributions
\[
\cA h=-\sum_{j=1}^3\partial_j(\omega M\partial_jh)+J^T\curl(\omega\br\bt\Omega_h).
\]
Furthermore,
\begin{equation}\label{eq:invnorm}
\cE_-(g):=\frac12\ip g{\cA^{-1}g}=\frac12a(h,h)\asymp\norm g_{\Hm}^2,
\qquad\norm{\nabla h}_2+\norm h_6\le C\norm g_{\Hm}.
\end{equation}
\end{lemma}

\begin{proof}
The sum-of-squares formula for $(\partial_jh)^TM\partial_jh$, uniform positivity of $(\br,\bt,\omega)$, and $|\Omega_h|\le C|\nabla h|$ imply
\[
c_A\norm{\nabla h}_2^2\le a(h,h)\le C_A\norm{\nabla h}_2^2
\qquad\text{for all }h\in\Hp(\R^3;\R^5).
\]
By the Riesz representation theorem on $\Hp(\R^3;\R^5)$ and the Sobolev embedding $\Hp(\R^3)\hookrightarrow L^6(\R^3)$, each $g\in\Hm(\R^3;\R^5)$ admits a unique solution $h\in\Hp(\R^3;\R^5)\subset L^6(\R^3;\R^5)$ to $a(h,\varphi)=\ip g\varphi$. Taking $\varphi=h$ gives
\[
c_A\norm{\nabla h}_2^2\le\ip gh\le\norm g_{\Hm}\norm{\nabla h}_2,
\qquad \ip gh\le c_A^{-1}\norm g_{\Hm}^2,
\]
while for any $\varphi\in\Hp(\R^3;\R^5)$ with $\norm{\nabla\varphi}_2=1$, Cauchy--Schwarz yields
\[
|\ip g\varphi|=|a(h,\varphi)|\le C_A^{1/2}(\ip gh)^{1/2}\norm{\nabla\varphi}_2,
\]
so taking the supremum over $\norm{\nabla\varphi}_2=1$ gives $\norm g_{\Hm}^2\le C_A\ip gh$ and proves \eqref{eq:invnorm}.
\end{proof}

Set
\[
D_0=\norm{D^2h}_2^2,\qquad D_s=\int\omega|\bu'|\,|\nabla h|^2\dd.
\]

\begin{lemma}\label{lem:elliptic}
If $g\in L^2(\R^3)\cap\Hm(\R^3)$, then $h=\cA^{-1}g$ satisfies $D^2h\in L^2(\R^3)$. Once $k>0$ is fixed, there exists $\delta_0>0$ such that for all $0<\delta\le\delta_0$,
\begin{equation}\label{eq:elliptic}
\norm{D^2h}_2\le C\bigl(\norm g_2+\norm{\bu'\nabla h}_2\bigr),\qquad
\norm g_2^2\le C(D_0+\delta^2D_s).
\end{equation}
\end{lemma}

\begin{proof}
The fifth component of $J^T\curl(\omega\br\bt\Omega_h)$ equals
\[
\bu\bigl\{\partial_2(\omega\br\bt\Omega_{h,3})-\partial_3(\omega\br\bt\Omega_{h,2})\bigr\},
\]
which involves only the transverse derivatives $\partial_2,\partial_3$. Since $\bu=\bu(x_1)$ is independent of $(x_2,x_3)$, $\bu$ commutes with $\partial_2,\partial_3$, so $\cA$ can be written in divergence form as
\[
\cA h=-\sum_{i,j=1}^3\partial_i\bigl(C^{ij}(x_1)\partial_jh\bigr),
\]
where the $5\times5$ matrices $C^{ij}(x_1)$ correspond to the coercive form \eqref{eq:form}.

For any coordinate direction $e_\ell$ and step size $\tau\ne0$, define the difference quotient
\[
D_\tau^\ell h(x)=\frac{h(x+\tau e_\ell)-h(x)}{\tau}\in L^2(\R^3)\cap\Hp(\R^3),\qquad
\norm{D_\tau^\ell h}_2\le\norm{\partial_\ell h}_2.
\]
Testing $a(h,\varphi)=\ip g\varphi$ with $\varphi=-D_{-\tau}^\ell D_\tau^\ell h\in\Hp(\R^3)$ and using discrete integration by parts gives
\[
c_A\norm{\nabla D_\tau^\ell h}_2^2
\le C\norm{\nabla C}_\infty\norm{\nabla h}_2\norm{\nabla D_\tau^\ell h}_2
+\norm g_2\norm{\nabla D_\tau^\ell h}_2,
\]
and hence
\[
\norm{\nabla D_\tau^\ell h}_2\le C\bigl(\norm g_2+\norm{\nabla C}_\infty\norm{\nabla h}_2\bigr)
\]
uniformly in $\tau\ne0$. Letting $\tau\to0$ yields $D^2h\in L^2(\R^3)$.

At the constant reference state, the Fourier symbol of $\cA_*$ is $|\xi|^2M_*$ augmented by the nonnegative rotational matrix $\rho_*\theta_*(|\xi|^2I-\xi\otimes\xi)$ on the momentum block, so Plancherel's theorem gives $\norm{D^2v}_2\le C\norm{\cA_*v}_2$. Applying this bound to
\[
\cA_*h=g+\sum_{i,j=1}^3\partial_i\bigl\{(C^{ij}-C_*^{ij})\partial_jh\bigr\}
\]
and using $\norm{C-C_*}_\infty\le C_k\delta$ and $|(C^{ij})'|\le C_k|\bu'|$ gives
\[
\norm{D^2h}_2
\le C\norm g_2+C_k\delta\norm{D^2h}_2+C_k\norm{\bu'\nabla h}_2.
\]
Choosing $\delta>0$ small enough that $C_k\delta\le1/2$ proves the first inequality in \eqref{eq:elliptic}. Conversely, expanding $\cA h$ directly and using
\[
\norm g_2\le C\bigl(\norm{D^2h}_2+\norm{\bu'\nabla h}_2\bigr),\qquad
\norm{\bu'\nabla h}_2^2
\le C\norm{\bu'}_\infty\int\omega|\bu'|\,|\nabla h|^2\dd
\le C\delta^2D_s
\]
gives the second inequality in \eqref{eq:elliptic}. \leavevmode{Finally, since $\cA$ is independent of $t$ and the conservative flux identity above ensures that $g\in H^1(0,T;\Hm(\R^3))$, differentiating $\cE_-(g)=\frac12\ip g{\cA^{-1}g}$ in time yields}
\begin{equation}\label{eq:invder}
\cE_-'=\ip{g_t}h.
\end{equation}
\end{proof}

\section[Energy estimates]{{Energy estimates}}\label{sec:compression}

Recall that $g=(\rr,n,\EE)^T$ with $n=\mm+d\nabla\rr$ and $h=\cA^{-1}g$. We now establish the main negative-order energy inequality.

\begin{proposition}\label{prop:negative_energy}
Let $d>0$ be given by \eqref{eq:d}. There exist a fixed weight exponent $k>0$ and constants $c,C>0$, independent of $\delta$ and $N=\norm Z_{H^3}$, such that for all sufficiently small $\delta$ and $N$,
\begin{equation}\label{eq:negative_estimate}
\cE_-' +c(D_0+D_s)
\le C\norm{\nabla\rr}_2^2
+CN\left\{D_0+\norm Z_2^2+\norm{\nabla(\psi,\zeta)}_2^2\right\}.
\end{equation}
\end{proposition}

\begin{proof}
\leavevmode{By Lemmas~\ref{lem:inverse} and~\ref{lem:elliptic}, the potential $h=\cA^{-1}g$ satisfies $h\in L^6(\R^3)$, $\nabla h\in H^1(\R^3)\subset L^2(\R^3)\cap L^6(\R^3)$, and $D^2h\in L^2(\R^3)$. Since $C_c^\infty(\R^3)$ is dense in $\Hp(\R^3)\cap\dot H^2(\R^3)$ under the norm $\norm{\nabla\cdot}_2+\norm{D^2\cdot}_2$, all integrations by parts in the pairings below are justified by standard density arguments.}

\medskip\noindent
\noindent\textbf{Step 1. Effective perturbation system and energy identity.}
Substituting $f=g-(0,d\nabla\rr,0)^T$ into \eqref{eq:conservative_perturbation} and using $\rr_t-s\partial_1\rr=-\diver n+d\Delta\rr$, we find that $g$ solves
\begin{equation}\label{eq:effective_perturbation}
g_t-s\partial_1g+\sum_{j=1}^3\partial_j(A_jg)
=\mathcal B g+\cN_g.
\end{equation}
Here the linearized velocity and temperature variations $\hat u,\hat\theta$ and the profile-gradient operator $\mathcal T\rr$ are defined by
\begin{align*}
\hat u&=\frac{n-\bu\rr e_1-d\nabla\rr}{\br},\\
\hat\theta&=\frac{\EE-\bu n_1+(\bu^2/2-c_v\bt)\rr+d\bu\partial_1\rr}{c_v\br},\\
\mathcal T\rr&=\partial_1(\bu\nabla\rr)+e_1\diver(\bu\nabla\rr)
-(\gamma-1)\nabla(\bu\partial_1\rr),
\end{align*}
and the linear diffusive operator $\mathcal Bg$ takes the explicit form
\begin{equation}\label{eq:rhs}
\mathcal B g=\begin{pmatrix}
d\Delta\rr\\
\cL\hat u-d\nabla\diver n+d^2\nabla\Delta\rr+d\mathcal T\rr\\
\kappa\Delta\hat\theta+\diver\{S(\bar u)\hat u+S(\hat u)\bar u\}
+d\diver(i_0\nabla\rr)-d(\gamma-1)\partial_1(\bu^2\partial_1\rr)
\end{pmatrix}.
\end{equation}
The exact quadratic remainders
\[
\mathcal Q_u=\psi-\hat u,\qquad
\mathcal Q_\theta=\zeta-\hat\theta,\qquad
\mathcal Q_{F,j}=F_j(\bar U+f)-F_j(\bar U)-A_jf=(0,\mathcal Q_{m,j},\mathcal Q_{E,j})^T
\]
are given by
\begin{align}
\mathcal Q_u&=-\frac{\rr}{\br}\psi,&
\mathcal Q_\theta&=-\frac{\rr}{\br}\zeta-\frac{\rho}{2c_v\br}|\psi|^2,\nonumber\\
\mathcal Q_m&=\rho\psi\otimes\psi-\frac{\rho}{2c_v}|\psi|^2I,&\label{eq:remainders}\\
\mathcal Q_E&=\rho\left\{(c_v+1)\zeta+\bar u\cdot\psi+\frac12|\psi|^2\right\}\psi
-\frac{\rho}{2c_v}|\psi|^2\bar u,&&\nonumber
\end{align}
and assemble into the divergence-form nonlinearity
\begin{equation}\label{eq:effective_nonlinearity}
\cN_g=\begin{pmatrix}
0\\
-\diver\mathcal Q_m+\cL\mathcal Q_u\\
-\diver\mathcal Q_E+\kappa\Delta\mathcal Q_\theta
+\diver\{S(\mathcal Q_u)\bar u+S(\bar u)\mathcal Q_u+S(\psi)\psi\}
\end{pmatrix}.
\end{equation}
Setting $\widetilde A_j=A_j-s\delta_{j1}I$ and $Y=\sum_{j=1}^3\widetilde A_j^T\partial_jh$, we pair \eqref{eq:effective_perturbation} with $h=\cA^{-1}g$ and use \eqref{eq:invder} to obtain
\begin{equation}\label{eq:negative_master}
\cE_-'=Q_E+V+Q_N,\qquad
Q_E=\ip{\cA h}Y,\quad V=\ip{\mathcal B g}h,\quad Q_N=\ip{\cN_g}h.
\end{equation}

\medskip\noindent
\noindent\textbf{Step 2. Coercivity of the Euler pairing $Q_E$.}
As shown in \eqref{eq:idea_neg_entropy_cancel} and \eqref{eq:idea_neg_curl_cancel}--\eqref{eq:idea_neg_curl_pairing}, the entropy symmetry $M\widetilde A_j^T=\widetilde A_jM$ and the curl identity $\curl\nabla=0$ eliminate all second derivatives of $h$ from $Q_E=Q_M+Q_\Omega$, yielding
\begin{equation}\label{eq:Qbase}
Q_M=-\frac12\int\sum_{i=1}^3
(\partial_ih)^T(\omega\widetilde A_1M)'\partial_ih\dd
+\int\omega\sum_{j=1}^3(\partial_1h)^TM(A_j')^T\partial_jh\dd,
\end{equation}
and
\begin{align}\label{eq:Qcurl}
Q_\Omega={}&-\frac12\int[(\bu-s)\omega\br\bt]'|\Omega_h|^2\dd\nonumber\\
&+\int\omega\br\bt\bu'\sum_{j=2}^3
\bigl\{(\partial_1h_{m,j})^2-(\partial_jh_{m,1}+\bu\partial_jh_E)^2\bigr\}\dd\\
&+\int\omega\br\bt(i_0'-s\bu')
\Omega_h\cdot(e_1\times\nabla h_E)\dd.\nonumber
\end{align}
\leavevmode{Using the change of variables preceding \eqref{eq:idea_neg_profile_comparison}, which rescales both Euler pairings by the same positive factor, and writing $c=\sqrt\gamma$ in the normalized variables, we see from \eqref{eq:idea_neg_phi}--\eqref{eq:idea_neg_square} that the weight-derivative terms at the reference state \eqref{eq:idea_neg_tangent} equal $\frac k2\omega\bu'F(\nabla h)$, where}
\begin{equation}\label{eq:F}
\begin{split}
F(\nabla h)=|\curl h_m|^2+\sum_{j=1}^3\Bigg\{
&\left(\frac{\partial_jh_\rho+(c_v+1)\partial_jh_E}c-\partial_jh_{m,1}\right)^2\\
&+\frac{(\partial_jh_\rho)^2}{c_v+1}+(\partial_jh_{m,2})^2+(\partial_jh_{m,3})^2
\Bigg\},
\end{split}
\end{equation}
whose null space is the one-dimensional normal acoustic line
\begin{equation}\label{eq:kernel}
\ker F=\{(\alpha\ell,0,0):\alpha\in\R\},\qquad
\ell=\bigl(0,\,1,\,0,\,0,\,(c_vc)^{-1}\bigr)^T.
\end{equation}
By \eqref{eq:idea_neg_metric_cancel}--\eqref{eq:idea_neg_coercivity}, the profile-derivative quadratic form $P(\nabla h)$ is strictly positive on $\ker F$, so fixing $k>0$ sufficiently large ensures
\[
P(\nabla h)+\frac k2F(\nabla h)\ge 2c_E|\nabla h|^2
\]
for a constant $c_E>0$. Since Lemma~\ref{lem:profile} and $i_0'-s\bu'=O(\delta|\bu'|)$ imply that the actual profile coefficients in \eqref{eq:Qbase}--\eqref{eq:Qcurl} differ from their values at the reference state by $O_k(\delta)$, choosing $\delta>0$ small enough that $C_k\delta\le c_E$ gives
\begin{equation}\label{eq:compression}
Q_E\le(-2c_E+C_k\delta)D_s\le-c_ED_s.
\end{equation}

\medskip\noindent
\noindent\textbf{Step 3. Constant-state parabolic coercivity.}\phantomsection\label{sec:diffusion}
In the original physical units, set
\[
\mu_*=\frac\mu{\rho_*},\qquad
a_L=\frac\nu{\rho_*}-d\ge\frac\nu{2\rho_*}>0,
\]
and extract from $\mathcal Bg$ the constant-state operator
\[
L_0g=\begin{pmatrix}
d\Delta\rr\\
\mu_*\Delta n+(a_L-\mu_*)\nabla\diver n\\
\theta_*[(c_v+1)d-\kappa/\rho_*]\Delta\rr
+\kappa(c_v\rho_*)^{-1}\Delta\EE
\end{pmatrix},\qquad
\ip{L_0g}h=\ip{\cA h}{L_0^*h}.
\]
On the scalar $(\rho,E)$ block, completing the square via Plancherel's identity gives the pairing with $\cA_*h$:
\begin{equation}\label{eq:scalar_square}
-\rho_*d\norm{D^2\bigl(h_\rho+(c_v+\tfrac12)\theta_*h_E\bigr)}_2^2
-\theta_*^2\left(\kappa-\frac14\rho_*d\right)\norm{D^2h_E}_2^2,
\end{equation}
where $\kappa-\frac14\rho_*d\ge\frac34\kappa>0$ by \eqref{eq:d}. On the momentum block, setting
\[
X_m=\nabla\diver h_m,\qquad C_m=\curl\curl h_m,
\]
we have
\[
(\cA_*h)_m=\rho_*\theta_*(-X_m+2C_m),\qquad
(L_0^*h)_m=a_LX_m-\mu_*C_m,\qquad
\ip{X_m}{C_m}=0.
\]
Using $\norm{D^2h_m}_2^2=\norm{X_m}_2^2+\norm{C_m}_2^2$, we obtain
\begin{equation}\label{eq:momentum_diff}
\ip{(\cA_*h)_m}{(L_0^*h)_m}
=-\rho_*\theta_*\bigl\{a_L\norm{X_m}_2^2+2\mu_*\norm{C_m}_2^2\bigr\},
\end{equation}
which combined with \eqref{eq:scalar_square} yields $\ip{\cA_*h}{L_0^*h}\le-c_0D_0$.

\medskip\noindent
\noindent\textbf{Step 4. Variable-coefficient diffusive and effective-momentum terms.}
By Lemma~\ref{lem:elliptic}, we have
\begin{gather*}
\norm{\bu'\nabla h}_2^2\le C\delta^2D_s,\qquad
\norm g_2\le C\bigl(D_0^{1/2}+\norm{\bu'\nabla h}_2\bigr),\\
\bigl|\ip{(\cA-\cA_*)h}{L_0^*h}\bigr|\le C\delta D_0+C\norm{\bu'\nabla h}_2D_0^{1/2}.
\end{gather*}
Decomposing $\hat u$ and $\hat\theta$ into their zeroth-order parts in $g$ and their density-gradient parts,
\[
\hat u=\hat u_0-\frac d{\br}\nabla\rr,\quad
\hat u_0=\frac{n-\bu\rr e_1}{\br},\qquad
\hat\theta=\hat\theta_0+\frac{d\bu}{c_v\br}\partial_1\rr,\quad
\hat\theta_0=\frac{\EE-\bu n_1+(\bu^2/2-c_v\bt)\rr}{c_v\br},
\]
and integrating two derivatives onto $h$ gives
\[
\bigl|\ip{\cL(\hat u_0-n/\rho_*)}{h_m}\bigr|
+\bigl|\ip{\kappa\Delta(\hat\theta_0-\EE/(c_v\rho_*)+\theta_*\rr/\rho_*)}{h_E}\bigr|
\le C\delta\norm g_2D_0^{1/2}.
\]
For the linearized viscous work in the energy equation, since $S(\bar u)=\diag(\nu\bu',\lambda\bu',\lambda\bu')$, integration by parts for any vector field $w$ cancels the normal derivative $\nu\bu'w_1\partial_1h_E$:
\begin{align}\label{eq:work}
\ip{\diver\{S(\bar u)w+S(w)\bar u\}}{h_E}
={}&\int\bu w\cdot\bigl\{(\mu+\lambda)\nabla\partial_1h_E
+\mu e_1\Delta h_E\bigr\}\dd\nonumber\\
&+(\mu-\lambda)\int\bu'w_\perp\cdot\nabla_\perp h_E\dd.
\end{align}
\leavevmode{For the $d$-dependent Euler terms in $\mathcal Bg$, setting $i_*=(c_v+1)\theta_*$ and integrating by parts twice to move the derivatives from $\rr$ onto $h_m$ and $h_E$ gives}
\begingroup\begin{align*}
\ip{\mathcal T\rr}{h_m}
&=\int\rr\bu\bigl\{\Delta h_{m,1}+(2-\gamma)\partial_1\diver h_m\bigr\}\dd\\
&\quad+\int\rr\bu'\bigl\{2\partial_1h_{m,1}-(\gamma-1)\diver h_m\bigr\}\dd,\\
\ip{\diver\{(i_0-i_*)\nabla\rr\}}{h_E}
&=\int\rr\bigl\{(i_0-i_*)\Delta h_E+i_0'\partial_1h_E\bigr\}\dd,\\
\ip{\partial_1(\bu^2\partial_1\rr)}{h_E}
&=\int\rr\bigl(\bu^2\partial_1^2h_E+2\bu\bu'\partial_1h_E\bigr)\dd.
\end{align*}
\endgroup
\leavevmode{Since $|\bu|+|i_0-i_*|\le C\delta$ and $|i_0'|\le C|\bu'|$, these three pairings are bounded by $C\norm g_2(\delta D_0^{1/2}+\norm{\bu'\nabla h}_2)$. Combining this estimate with \eqref{eq:work} evaluated at $w=\hat u_0$ and $w=-d\br^{-1}\nabla\rr$ bounds all variable-coefficient remainders considered above by}
\[
C\bigl(\delta\norm g_2+\norm{\nabla\rr}_2\bigr)\bigl(D_0^{1/2}+\norm{\bu'\nabla h}_2\bigr)+C\norm g_2\norm{\bu'\nabla h}_2.
\]
Moreover, the third-order density terms in \eqref{eq:rhs} integrate by parts onto $D^2h_m$ and $D^2h_E$:
\begin{equation}\label{eq:density3}
\begin{split}
\ip{-d\cL(\br^{-1}\nabla\rr)+d^2\nabla\Delta\rr}{h_m}
&=-d\int\br^{-1}\nabla\rr\cdot\cL h_m\dd
+d^2\int\nabla\rr\cdot\Delta h_m\dd,\\
\bigl|\ip{-d\cL(\br^{-1}\nabla\rr)+d^2\nabla\Delta\rr}{h_m}\bigr|
&\le C\norm{\nabla\rr}_2D_0^{1/2},
\end{split}
\end{equation}
and
\[
\left|\ip{\kappa d\Delta\left(\frac{\bu}{c_v\br}\partial_1\rr\right)}{h_E}\right|
\le C\delta\norm{\nabla\rr}_2D_0^{1/2}.
\]
\leavevmode{Combining the estimates in Steps 3 and 4 with $\norm g_2\le C(D_0^{1/2}+\norm{\bu'\nabla h}_2)$ and $\norm{\bu'\nabla h}_2^2\le C\delta^2D_s$, we find that}
\begingroup\[
V\le-(c_0-C\delta)D_0+CD_0^{1/2}\bigl(\delta D_s^{1/2}+\norm{\nabla\rr}_2\bigr)+C\delta^2D_s.
\]
\endgroup
\leavevmode{Applying Young's inequality to absorb the $D_0^{1/2}$ product terms into $-\frac12c_0D_0$ and taking $\delta>0$ sufficiently small yields}
\begin{equation}\label{eq:diffusion}
V\le-cD_0+\frac{c_E}{2}D_s+C\norm{\nabla\rr}_2^2.
\end{equation}

\medskip\noindent
\noindent\textbf{Step 5. Control of the nonlinear flux pairing $Q_N$.}
By \eqref{eq:remainders}, we have
\[
|\mathcal Q_u|+|\mathcal Q_\theta|+|\mathcal Q_F|\le C|f|^2,\qquad
\norm f_3+\norm\psi_3\le CN,\qquad
\norm{\nabla h}_6\le CD_0^{1/2}.
\]
\leavevmode{Moving the outer flux derivatives onto $h$, using \eqref{eq:work} with $w=\mathcal Q_u$, and retaining the first derivative in $S(\psi)$ gives}
\begin{align*}
|\ip{-\diver\mathcal Q_F}h|&\le\norm{\mathcal Q_F}_{6/5}\norm{\nabla h}_6\le CN\norm f_2D_0^{1/2},\\
|\ip{\cL\mathcal Q_u}{h_m}|+|\ip{\kappa\Delta\mathcal Q_\theta}{h_E}|&\le C(\norm{\mathcal Q_u}_2+\norm{\mathcal Q_\theta}_2)D_0^{1/2}\le CN\norm{\nabla f}_2D_0^{1/2},\\
|\ip{\diver(S(\psi)\psi)}{h_E}|&\le C\norm{\nabla\psi}_2\norm\psi_3\norm{\nabla h_E}_6\le CN\norm{\nabla\psi}_2D_0^{1/2},\\
\bigl|\ip{\diver\{S(\mathcal Q_u)\bar u+S(\bar u)\mathcal Q_u\}}{h_E}\bigr|&\le CN\bigl(\delta\norm{\nabla f}_2+\delta^2\norm f_2\bigr)D_0^{1/2},
\end{align*}
whence
\begin{equation}\label{eq:nonlinear}
|Q_N|=|\ip{\cN_g}h|
\le CND_0^{1/2}\bigl(\norm f_2+\norm{\nabla f}_2+\norm{\nabla\psi}_2\bigr).
\end{equation}
Since
\[
\norm f_2\le C\norm Z_2,\qquad
\norm{\nabla f}_2\le C\bigl(\norm{\nabla\rr}_2+\norm{\nabla(\psi,\zeta)}_2+\delta^2\norm Z_2\bigr),
\]
combining \eqref{eq:compression}, \eqref{eq:diffusion}, and \eqref{eq:nonlinear} in \eqref{eq:negative_master} proves \eqref{eq:negative_estimate}.
\end{proof}

\begin{lemma}\label{lem:negative_recovery}
For $m=1,2,3$, define the dissipation
\[
D_m(Z)=\norm{\nabla\phi}_{H^{m-1}}^2+\norm{\nabla(\psi,\zeta)}_{H^m}^2.
\]
Under \eqref{eq:bootstrap}, one has
\begin{align}
\norm f_2^2+\norm Z_{H^m}^2
&\le C\bigl\{D_0+D_m(Z)+\delta^2D_s\bigr\},\label{eq:recovery}\\
\norm{\nabla\rr}_2^2
&\le CD_m(Z)+C\delta^4(D_0+\delta^2D_s),\label{eq:rho_recovery}\\
\norm{\nabla f}_2^2
&\le CD_m(Z)+C\delta^4\norm Z_2^2.\nonumber
\end{align}
\end{lemma}

\begin{proof}
From the exact relations
\begin{gather*}
g_\rho=\rr=\br(e^\phi-1),\qquad
\nabla\rr=\rho\nabla\phi+\br'(e^\phi-1)e_1,\qquad
\rho\psi=n-\bar u\rr-d\nabla\rr,\\
c_v\rho\zeta=\EE-\rr\left(c_v\bt+\frac12|\bar u|^2\right)-\rho\bar u\cdot\psi-\frac12\rho|\psi|^2,
\end{gather*}
we obtain
\begin{gather*}
\norm\phi_2\le C\norm\rr_2\le C\norm g_2,\qquad
\norm{\nabla\rr}_2\le C\norm{\nabla\phi}_2+C\delta^2\norm g_2,\\
\norm\psi_2+\norm\zeta_2\le C\bigl(\norm g_2+\norm{\nabla\phi}_2\bigr).
\end{gather*}
Combining these bounds with the elliptic estimate \eqref{eq:elliptic} and $\norm Z_{H^m}^2\le C(\norm Z_2^2+D_m(Z))$ yields \eqref{eq:recovery}--\eqref{eq:rho_recovery}. {Finally, differentiating the algebraic expressions for $f=(\rr,\mm,\EE)^T$ in terms of $(\rr,\psi,\zeta)$ and $\bar U$, and using $|\bar U'|\le C|\bu'|\le C\delta^2$ together with $\norm{\nabla\rr}_2\le C(\norm{\nabla\phi}_2+\delta^2\norm Z_2)$, we obtain the last bound for $\norm{\nabla f}_2^2$.}
\end{proof}

We now establish the energy estimates for $Z=(\phi,\psi,\zeta)^T$ of orders $m=0,1,2,3$ under the a priori assumption \eqref{eq:bootstrap}.

\leavevmode{From now on, we set}
\[
a=\rho^{-1},\qquad \widetilde a=a-\br^{-1},\qquad
D_t=\partial_t+(u-se_1)\cdot\nabla,\qquad N=\norm Z_{H^3}.
\]
For any multi-index $\alpha\in\mathbb N_0^3$, we write $z_\alpha=\partial^\alpha z$ and denote the commutator by $[\partial^\alpha,b]w=\partial^\alpha(bw)-b\partial^\alpha w$. For any perturbation triple $V=(\phi,\psi,\zeta)^T$, define
\begin{align}\label{eq:Pm}
S_m(V)&=\frac12\sum_{|\alpha|\le m}\int
\left(\theta_*\phi_\alpha^2+|\psi_\alpha|^2+
\frac{c_v}{\theta_*}\zeta_\alpha^2\right)\dd,
&&0\le m\le3,\nonumber\\
I_m(V)&=\sum_{|\alpha|\le m-1}\int \psi_\alpha\cdot\nabla\phi_\alpha\dd,
\qquad P_m(V)=S_m(V)+\sigma_m I_m(V),&&1\le m\le3,\\
D_m(V)&=\norm{\nabla\phi}_{H^{m-1}}^2+
\norm{\nabla(\psi,\zeta)}_{H^m}^2,&&1\le m\le3.\nonumber
\end{align}
By Lemma~\ref{lem:profile} and Sobolev embedding under \eqref{eq:bootstrap}, the coefficients satisfy
\begin{gather}
\norm{\theta-\theta_*}_\infty\le C(\delta+N),\qquad
\norm{\nabla(u,\theta,a)}_\infty\le C(\delta^2+N),\nonumber\\
\norm{\widetilde a}_{H^m}\le C\norm Z_{H^m}\quad(0\le m\le3),\qquad
0<c\le\rho,\theta,a\le C,\nonumber\\
|I_m(Z)|\le C\norm Z_{H^m}^2,\qquad
\norm Z_{H^m}^2\le C\bigl(D_m(Z)+\norm Z_2^2\bigr)\quad(1\le m\le3).
\label{eq:physical_coeff_bounds}
\end{gather}

Subtracting the stationary shock profile from \eqref{eq:nsf} gives the primitive perturbation system
\begin{equation}\label{eq:physical}
\left\{
\begin{aligned}
&D_t\phi+\diver\psi=F_\phi,\\
&D_t\psi+\theta\nabla\phi+\nabla\zeta-a\cL\psi=F_\psi,\\
&c_vD_t\zeta+\theta\diver\psi-\kappa a\Delta\zeta=F_\zeta,
\end{aligned}
\right.
\end{equation}
with source terms
\begin{align}\label{eq:physical_sources}
F_\phi&=-\frac{\br'}{\br}\psi_1,\nonumber\\
F_\psi&=-\bu'\psi_1e_1-\frac{\br'}{\br}\zeta e_1+\nu\bu''\widetilde a e_1,\\
F_{\zeta,0}&=-c_v\bt'\psi_1-\bu'\zeta+\kappa\widetilde a\bt''+\nu(\bu')^2\widetilde a,\nonumber\\
F_\zeta&=F_{\zeta,0}+\mathcal H,\qquad
\mathcal H=aS(\psi):\nabla\psi+2aS(\bar u):\nabla\psi.\nonumber
\end{align}
Applying $\partial^\alpha$ ($|\alpha|\le 3$) to \eqref{eq:physical} yields the differentiated system
\begin{equation}\label{eq:physical_differentiated}
\left\{
\begin{aligned}
&D_t\phi_\alpha+\diver\psi_\alpha=F_{\phi,\alpha},\\
&D_t\psi_\alpha+\theta\nabla\phi_\alpha+\nabla\zeta_\alpha-a\cL\psi_\alpha=F_{\psi,\alpha},\\
&c_vD_t\zeta_\alpha+\theta\diver\psi_\alpha-\kappa a\Delta\zeta_\alpha=F_{\zeta,\alpha},
\end{aligned}
\right.
\end{equation}
where the commutator and source terms are
\begin{align}\label{eq:physical_commutators}
F_{\phi,\alpha}&=\partial^\alpha F_\phi-[\partial^\alpha,u]\cdot\nabla\phi,\nonumber\\
F_{\psi,\alpha}&=\partial^\alpha F_\psi-[\partial^\alpha,u]\cdot\nabla\psi
-[\partial^\alpha,\theta]\nabla\phi+[\partial^\alpha,a]\cL\psi,\\
F_{\zeta,\alpha}&=\partial^\alpha F_\zeta-c_v[\partial^\alpha,u]\cdot\nabla\zeta
-[\partial^\alpha,\theta]\diver\psi
+\kappa[\partial^\alpha,a]\Delta\zeta.\nonumber
\end{align}

Testing \eqref{eq:physical_differentiated} with the constant-state multipliers
\begin{equation}\label{eq:l2_multipliers}
\theta_*\phi_\alpha,\qquad\psi_\alpha,\qquad\frac{\zeta_\alpha}{\theta_*}
\end{equation}
and summing over $|\alpha|\le m$ ($0\le m\le 3$), we see that the constant-state acoustic couplings vanish for each $\alpha$:
\begin{equation}\label{eq:physical_acoustic_cancel}
\int\theta_*(\phi_\alpha\diver\psi_\alpha+\psi_\alpha\cdot\nabla\phi_\alpha)\dd=0,
\qquad
\int(\psi_\alpha\cdot\nabla\zeta_\alpha+\zeta_\alpha\diver\psi_\alpha)\dd=0.
\end{equation}
With the parabolic dissipation
\[
\mathcal V_m=\sum_{|\alpha|\le m}\int a
\left\{\mu|\nabla\psi_\alpha|^2+(\mu+\lambda)|\diver\psi_\alpha|^2+
\frac\kappa{\theta_*}|\nabla\zeta_\alpha|^2\right\}\dd
\ge c_0\norm{\nabla(\psi,\zeta)}_{H^m}^2,
\]
integration by parts over $\R^3$ yields the symmetric energy identity
\begin{equation}\label{eq:fixed_symmetric_identity}
S_m(Z)'+\mathcal V_m=\sum_{|\alpha|\le m}\bigl(\mathcal R_{1,\alpha}+\mathcal R_{2,\alpha}+\mathcal R_{3,\alpha}+\mathcal R_{4,\alpha}\bigr),
\end{equation}
which at $m=0$ reads
\begin{equation}\label{eq:physical_l2_identity}
S_0(Z)'+\mathcal V_0=\mathcal R_{1,0}+\mathcal R_{2,0}+\mathcal R_{3,0}+\mathcal R_{4,0},
\end{equation}
where for each $\alpha$ the four remainder integrals are
\begin{align*}
\mathcal R_{1,\alpha}
&=\frac12\int(\diver u)
\left(\theta_*\phi_\alpha^2+|\psi_\alpha|^2+
\frac{c_v}{\theta_*}\zeta_\alpha^2\right)\dd,\\
\mathcal R_{2,\alpha}
&=-\int(\theta-\theta_*)\psi_\alpha\cdot\nabla\phi_\alpha\dd
-\int\frac{\theta-\theta_*}{\theta_*}\zeta_\alpha\diver\psi_\alpha\dd,\\
\mathcal R_{3,\alpha}
&=-\int\left\{\mu\sum_{j=1}^3(\partial_ja)\psi_\alpha\cdot\partial_j\psi_\alpha
+(\mu+\lambda)(\nabla a\cdot\psi_\alpha)\diver\psi_\alpha
+\frac\kappa{\theta_*}\zeta_\alpha\nabla a\cdot\nabla\zeta_\alpha\right\}\dd,\\
\mathcal R_{4,\alpha}
&=\int\left(\theta_*\phi_\alpha F_{\phi,\alpha}+
\psi_\alpha\cdot F_{\psi,\alpha}+\frac{\zeta_\alpha}{\theta_*}F_{\zeta,\alpha}\right)\dd.
\end{align*}
In $\mathcal R_{2,\alpha}$, we integrate the pressure term by parts so that no extra derivative falls on $\phi_\alpha$:
\begin{equation}\label{eq:pressure_transfer}
-\int(\theta-\theta_*)\psi_\alpha\cdot\nabla\phi_\alpha\dd
=\int(\theta-\theta_*)\phi_\alpha\diver\psi_\alpha\dd
+\int\phi_\alpha\psi_\alpha\cdot\nabla\theta\dd.
\end{equation}

For the cross functional $I_m(Z)$ ($1\le m\le 3$), pairing the momentum equation in \eqref{eq:physical_differentiated} with $\nabla\phi_\alpha$ and the gradient of the continuity equation with $\psi_\alpha$ for $|\alpha|\le m-1$ gives
\begin{equation}\label{eq:fixed_compensator_identity}
I_m(Z)'+\sum_{|\alpha|\le m-1}\int\theta|\nabla\phi_\alpha|^2\dd
=\sum_{|\alpha|\le m-1}\mathcal K_\alpha,
\end{equation}
where
\begin{align*}
\mathcal K_\alpha={}&-\int\nabla\zeta_\alpha\cdot\nabla\phi_\alpha\dd
+\int a\cL\psi_\alpha\cdot\nabla\phi_\alpha\dd
+\norm{\diver\psi_\alpha}_2^2\\
&+\int\left\{(\diver u)\psi_\alpha\cdot\nabla\phi_\alpha
-\psi_\alpha\cdot(\nabla u)^T\nabla\phi_\alpha
+F_{\psi,\alpha}\cdot\nabla\phi_\alpha-(\diver\psi_\alpha)F_{\phi,\alpha}\right\}\dd.
\end{align*}

\begin{lemma}\label{lem:physical_l2}
For sufficiently small $\delta+N$, the zeroth-order physical energy satisfies
\begin{equation}\label{eq:physical_l2_estimate}
S_0(Z)'+c\norm{\nabla(\psi,\zeta)}_2^2\le C(\delta+N)\norm Z_2^2.
\end{equation}
\end{lemma}

\begin{proof}
In \eqref{eq:physical_l2_identity}, using \eqref{eq:physical_coeff_bounds} and \eqref{eq:pressure_transfer} at $\alpha=0$, we obtain
\begin{align*}
|\mathcal R_{1,0}|&\le C(\delta^2+N)\norm Z_2^2,\\
|\mathcal R_{2,0}|+|\mathcal R_{3,0}|&\le C(\delta+N)\norm Z_2\norm{\nabla(\psi,\zeta)}_2+C(\delta^2+N)\norm Z_2^2.
\end{align*}
By Lemma~\ref{lem:profile} and \eqref{eq:physical_sources}, the source terms satisfy
\[
\norm{(F_\phi,F_\psi,F_{\zeta,0})}_2\le C\delta^2\norm Z_2,\qquad
\left|\int\frac{\zeta}{\theta_*}\mathcal H\dd\right|\le CN\norm{\nabla\psi}_2^2+C\delta^2\norm Z_2\norm{\nabla\psi}_2.
\]
Combining these bounds in \eqref{eq:physical_l2_identity} and absorbing $C(\delta+N)\norm{\nabla(\psi,\zeta)}_2^2$ into $\mathcal V_0\ge c_0\norm{\nabla(\psi,\zeta)}_2^2$ gives \eqref{eq:physical_l2_estimate}.
\end{proof}

\begin{lemma}\label{lem:physical_h1}
There exists a fixed constant $\sigma_1>0$ such that $P_1(Z)\asymp\norm Z_{H^1}^2$ and, for sufficiently small $\delta+N$,
\begin{equation}\label{eq:h1_physical_energy}
P_1(Z)'+c\bigl\{\norm{\nabla\phi}_2^2+\norm{\nabla(\psi,\zeta)}_{H^1}^2\bigr\}
\le C(\delta+N)\norm Z_2^2.
\end{equation}
\end{lemma}

\begin{proof}
In \eqref{eq:fixed_symmetric_identity} with $m=1$, \eqref{eq:physical_coeff_bounds} and \eqref{eq:pressure_transfer} yield
\[
\sum_{|\alpha|\le1}\bigl(|\mathcal R_{1,\alpha}|+|\mathcal R_{2,\alpha}|+|\mathcal R_{3,\alpha}|\bigr)
\le C(\delta+N)\bigl(\norm Z_{H^1}^2+\norm{\nabla(\psi,\zeta)}_{H^1}^2\bigr).
\]
Using $\norm{\mathcal H}_2\le C(\delta^2+N)\norm Z_{H^1}$, we bound the undifferentiated sources by
\begin{equation}\label{eq:first_order_physical_sources}
\norm{F_\phi}_{H^1}+\norm{F_\psi}_2+\norm{F_\zeta}_2\le C(\delta^2+N)\norm Z_{H^1}.
\end{equation}
Integrating by parts so that no derivative acts on $F_\psi$ or $F_\zeta$, we have
\[
\sum_{|\alpha|\le1}\ip{\partial^\alpha F_\psi}{\psi_\alpha}=\ip{F_\psi}{\psi-\Delta\psi},\qquad
\sum_{|\alpha|\le1}\ip{\partial^\alpha F_\zeta}{\zeta_\alpha}=\ip{F_\zeta}{\zeta-\Delta\zeta},
\]
while for $|\alpha|=1$ the commutators in \eqref{eq:physical_commutators} satisfy
\begin{align*}
\norm{[\partial_j,u]\cdot\nabla Z}_2
+\norm{[\partial_j,\theta](\nabla\phi,\diver\psi)}_2&\le C(\delta^2+N)\norm Z_{H^1},\\
\norm{[\partial_j,a](\cL\psi,\kappa\Delta\zeta)}_2&\le C(\delta^2+N)\norm{\nabla(\psi,\zeta)}_{H^1}.
\end{align*}
Consequently,
\[
\left|\sum_{|\alpha|\le1}\mathcal R_{4,\alpha}\right|
\le C(\delta^2+N)\bigl(\norm Z_{H^1}^2+\norm{\nabla(\psi,\zeta)}_{H^1}^2\bigr),
\]
and \eqref{eq:fixed_symmetric_identity} gives
\begin{equation}\label{eq:first_symmetric_estimate}
S_1(Z)'+c_0\norm{\nabla(\psi,\zeta)}_{H^1}^2
\le C(\delta+N)\bigl(\norm Z_{H^1}^2+\norm{\nabla(\psi,\zeta)}_{H^1}^2\bigr).
\end{equation}

In the cross identity \eqref{eq:fixed_compensator_identity} with $m=1$, positivity of $\theta$ gives $\int\theta|\nabla\phi|^2\dd\ge c\norm{\nabla\phi}_2^2$, and Young's inequality together with \eqref{eq:first_order_physical_sources} yields
\begin{equation}\label{eq:fixed_compensator_estimate}
I_1(Z)'+c_1\norm{\nabla\phi}_2^2
\le C_1\norm{\nabla(\psi,\zeta)}_{H^1}^2+C(\delta^2+N)\bigl(\norm Z_{H^1}^2+\norm{\nabla(\psi,\zeta)}_{H^1}^2\bigr).
\end{equation}
Fixing $\sigma_1>0$ small enough that $P_1(Z)\asymp\norm Z_{H^1}^2$ and $\sigma_1C_1\le c_0/2$, adding $\sigma_1$ times \eqref{eq:fixed_compensator_estimate} to \eqref{eq:first_symmetric_estimate}, and using $\norm Z_{H^1}^2\le C(D_1(Z)+\norm Z_2^2)$ proves \eqref{eq:h1_physical_energy}.
\end{proof}

\begin{lemma}\label{lem:physical_h2}
There exists a fixed constant $\sigma_2>0$ such that $P_2(Z)\asymp\norm Z_{H^2}^2$ and, for sufficiently small $\delta+N$,
\begin{equation}\label{eq:h2_physical_energy}
P_2(Z)'+c\bigl\{\norm{\nabla\phi}_{H^1}^2+\norm{\nabla(\psi,\zeta)}_{H^2}^2\bigr\}
\le C(\delta+N)\norm Z_2^2.
\end{equation}
\end{lemma}

\begin{proof}
At order $m=2$, combining \eqref{eq:fixed_symmetric_identity} and \eqref{eq:fixed_compensator_identity} gives
\begin{equation}\label{eq:h2_energy_identity}
P_2(Z)'+\mathcal V_2+
\sigma_2\sum_{|\alpha|\le1}\int\theta|\nabla\phi_\alpha|^2\dd
=\sum_{|\alpha|\le2}\sum_{k=1}^4\mathcal R_{k,\alpha}
+\sigma_2\sum_{|\alpha|\le1}\mathcal K_\alpha.
\end{equation}
By \eqref{eq:physical_coeff_bounds} and \eqref{eq:pressure_transfer},
\[
\sum_{|\alpha|\le2}\bigl(|\mathcal R_{1,\alpha}|+|\mathcal R_{2,\alpha}|+|\mathcal R_{3,\alpha}|\bigr)
\le C(\delta+N)\bigl(\norm Z_{H^2}^2+\norm{\nabla(\psi,\zeta)}_{H^2}^2\bigr).
\]
In the second-order commutators, any quadratic product of perturbation derivatives $y,z\in\{\phi,\psi,\zeta\}$ satisfies
\[
\norm{\nabla y\,\partial^2z}_2+\norm{\partial^2y\,\nabla z}_2
\le\norm{\nabla y}_\infty\norm{\partial^2z}_2+\norm{\partial^2y}_2\norm{\nabla z}_\infty
\le CN\norm Z_{H^2},
\]
while for the diffusion commutators with $z\in\{\psi,\zeta\}$, writing $a=\br^{-1}+\widetilde a$ and using H\"older's inequality and Sobolev embeddings, we have
\[
\norm{\partial^2\widetilde a}_3\le CN,\qquad
\norm{\partial^2z}_6\le C\norm{\partial^3z}_2\le C\norm{\nabla(\psi,\zeta)}_{H^2},
\]
{where $\partial^k$ stands for any spatial derivative of order $k$. Since $\norm{\partial_1^j(\br^{-1})}_\infty\le C\delta^2$ ($j=1,2$) and $\norm{\nabla a}_\infty\le C(\delta^2+N)$, estimating $(\partial^2\widetilde a)\partial^2z$ in $L^3\times L^6$ and both $(\partial_1^2\br^{-1})\partial^2z$ and $(\nabla a)\partial^3z$ in $L^\infty\times L^2$ yields}
\[
\sum_{1\le|\alpha|\le2}\bigl|\ip{[\partial^\alpha,a]\partial^2z}{z_\alpha}\bigr|
\le C(\delta^2+N)\bigl(\norm Z_{H^2}^2+\norm{\nabla(\psi,\zeta)}_{H^2}^2\bigr).
\]
Furthermore, the source terms satisfy
\[
\norm{(F_\phi,F_\psi,F_{\zeta,0})}_{H^2}\le C\delta^2\norm Z_{H^2},\qquad
\norm{\mathcal H}_{H^2}\le C(\delta^2+N)\bigl(\norm{\nabla(\psi,\zeta)}_{H^2}+\norm Z_{H^2}\bigr),
\]
so that
\[
S_2(Z)'+c_0\norm{\nabla(\psi,\zeta)}_{H^2}^2
\le C(\delta+N)\bigl(\norm Z_{H^2}^2+\norm{\nabla(\psi,\zeta)}_{H^2}^2\bigr).
\]
In the cross functional $I_2(Z)$, where $|\alpha|\le1$, we have
\[
\sum_{|\alpha|\le1}\norm{(F_{\phi,\alpha},F_{\psi,\alpha})}_2
\le C(\delta^2+N)\bigl(\norm Z_{H^2}+\norm{\nabla(\psi,\zeta)}_{H^2}\bigr),
\]
and Young's inequality yields
\[
I_2(Z)'+c_1\norm{\nabla\phi}_{H^1}^2
\le C_1\norm{\nabla(\psi,\zeta)}_{H^2}^2+C(\delta^2+N)\bigl(\norm Z_{H^2}^2+\norm{\nabla(\psi,\zeta)}_{H^2}^2\bigr).
\]
Choosing $\sigma_2>0$ small enough that $P_2(Z)\asymp\norm Z_{H^2}^2$ and $\sigma_2C_1\le c_0/2$, and using $\norm Z_{H^2}^2\le C(D_2(Z)+\norm Z_2^2)$, we absorb $C(\delta+N)D_2(Z)$ to obtain \eqref{eq:h2_physical_energy}.
\end{proof}

\begin{lemma}\label{lem:physical_h3}
There exists a fixed constant $\sigma_3>0$ such that $P_3(Z)\asymp N^2$ and, for sufficiently small $\delta+N$,
\begin{equation}\label{eq:h3_physical_energy}
P_3(Z)'+c\bigl\{\norm{\nabla\phi}_{H^2}^2+\norm{\nabla(\psi,\zeta)}_{H^3}^2\bigr\}
\le C(\delta+N)\norm Z_2^2.
\end{equation}
\end{lemma}

\begin{proof}
At top order $m=3$, combining \eqref{eq:fixed_symmetric_identity} and \eqref{eq:fixed_compensator_identity} gives
\begin{equation}\label{eq:h3_energy_identity}
P_3(Z)'+\mathcal V_3+
\sigma_3\sum_{|\alpha|\le2}\int\theta|\nabla\phi_\alpha|^2\dd
=\sum_{|\alpha|\le3}\sum_{k=1}^4\mathcal R_{k,\alpha}
+\sigma_3\sum_{|\alpha|\le2}\mathcal K_\alpha.
\end{equation}
\leavevmode{For $|\alpha|=3$, the two terms on the right of \eqref{eq:pressure_transfer} are bounded in absolute value by}
\begingroup\[
\norm{\phi_\alpha}_2\big(\norm{\theta-\theta_*}_\infty\norm{\diver\psi_\alpha}_2+
\norm{\nabla\theta}_\infty\norm{\psi_\alpha}_2\big).
\]
\endgroup
\leavevmode{By \eqref{eq:physical_coeff_bounds} this gives}
\[
\sum_{|\alpha|\le3}\bigl(|\mathcal R_{1,\alpha}|+|\mathcal R_{2,\alpha}|+|\mathcal R_{3,\alpha}|\bigr)
\le C(\delta+N)\bigl(N^2+\norm{\nabla(\psi,\zeta)}_{H^3}^2\bigr).
\]
For the third-order transport and pressure commutators with $y,z\in\{\phi,\psi,\zeta\}$, H\"older's inequality and Sobolev embeddings give
\[
\norm{\nabla y}_\infty\norm{\partial^3z}_2
+\norm{\partial^2y}_3\norm{\partial^2z}_6
+\norm{\partial^3y}_2\norm{\nabla z}_\infty
\le C\norm y_{H^3}\norm z_{H^3}\le CN^2.
\]
For the diffusion commutator pairings with $z\in\{\psi,\zeta\}$ at $|\alpha|=3$, let $I_{3,j}(z)$ denote the sum of absolute pairings where $j$ derivatives ($j=1,2,3$) act on $a=\br^{-1}+\widetilde a$. Using
\[
\norm{\partial^3z}_6\le C\norm{\partial^4z}_2\le C\norm{\nabla(\psi,\zeta)}_{H^3},\qquad
\norm{\partial^2z}_3\le CN,\qquad
\norm{\widetilde a}_{H^3}\le CN,
\]
we obtain
\begin{align}\label{eq:fixed_diffusion_commutators}
I_{3,1}(z)&\le C\norm{\nabla a}_\infty\norm{\partial^4z}_2\norm{\partial^3z}_2
\le C(\delta^2+N)N\norm{\nabla(\psi,\zeta)}_{H^3},\nonumber\\
I_{3,2}(z)&\le C\norm{\partial^2\widetilde a}_3\norm{\partial^3z}_6\norm{\partial^3z}_2+C\delta^2N^2
\le CN^2\norm{\nabla(\psi,\zeta)}_{H^3}+C\delta^2N^2,\\
I_{3,3}(z)&\le C\norm{\partial^3\widetilde a}_2\norm{\partial^2z}_3\norm{\partial^3z}_6+C\delta^2N^2
\le CN^2\norm{\nabla(\psi,\zeta)}_{H^3}+C\delta^2N^2.\nonumber
\end{align}
Moreover,
\begingroup\[
\norm{(F_\phi,F_\psi,F_{\zeta,0})}_{H^3}\le C\delta^2N,\qquad
\norm{\mathcal H}_{H^3}\le C(\delta^2+N)\norm{\nabla\psi}_{H^3}.
\]
\endgroup
\leavevmode{Indeed, decomposing $a=\br^{-1}+\widetilde a$ in \eqref{eq:physical_sources} and using $\norm{\br^{-1}}_{W^{3,\infty}}\le C$, $\norm{S(\bar u)}_{W^{3,\infty}}\le C\delta^2$, $\norm{\widetilde a}_{H^3}\le CN$, and $\norm{\nabla\psi}_\infty\le C\norm{\nabla\psi}_{H^2}\le CN$ (with $\norm{\nabla\psi}_\infty\le C\norm{\nabla\psi}_{H^3}$ when all three derivatives act on $\widetilde a$), standard Moser-type product estimates give the stated bound on $\norm{\mathcal H}_{H^3}$, whose $H^3$ pairing with $\zeta$ is therefore bounded by $C(\delta^2+N)N\norm{\nabla\psi}_{H^3}$. Combining these estimates yields}
\[
S_3(Z)'+c_0\norm{\nabla(\psi,\zeta)}_{H^3}^2
\le C(\delta+N)\bigl(N^2+\norm{\nabla(\psi,\zeta)}_{H^3}^2\bigr).
\]
In the cross functional $I_3(Z)$, where $|\alpha|\le2$, we have
\[
\sum_{|\alpha|\le2}\norm{(F_{\phi,\alpha},F_{\psi,\alpha})}_2
\le C(\delta^2+N)\bigl(N+\norm{\nabla(\psi,\zeta)}_{H^3}\bigr),
\]
and hence
\[
I_3(Z)'+c_1\norm{\nabla\phi}_{H^2}^2
\le C_1\norm{\nabla(\psi,\zeta)}_{H^3}^2+C(\delta^2+N)\bigl(N^2+\norm{\nabla(\psi,\zeta)}_{H^3}^2\bigr).
\]
Fixing $\sigma_3>0$ small enough that $P_3(Z)\asymp N^2$ and $\sigma_3C_1\le c_0/2$, and using $N^2\le C(D_3(Z)+\norm Z_2^2)$, we absorb $C(\delta+N)D_3(Z)$ to establish \eqref{eq:h3_physical_energy}.
\end{proof}

\begin{proposition}\label{prop:physical}
There exist uniform constants $c_p,c_0,C_p,C>0$ such that for sufficiently small $\delta+N$ and each $m=1,2,3$,
\begin{align}
P_m(Z)'+c_pD_m(Z)&\le C_p(\delta+N)\bigl\{D_m(Z)+\norm Z_2^2\bigr\},\label{eq:phys_est}\\
S_m(Z)'+c_0\norm{\nabla(\psi,\zeta)}_{H^m}^2&\le C(\delta+N)\bigl\{\norm Z_{H^m}^2+\norm{\nabla(\psi,\zeta)}_{H^m}^2\bigr\}.\label{eq:fixed_symmetric_estimate}
\end{align}
\end{proposition}

\leavevmode{For the time-decay analysis below, we record the following linear first-order estimate.}

\begin{lemma}\label{lem:general_physical}
\leavevmode{On a finite interval, let $V=(\dot\phi,\dot\psi,\dot\zeta)^T$ be a smooth solution of the linear system}
\[
\left\{
\begin{aligned}
&D_t\dot\phi+\diver\dot\psi=R_\phi,\\
&D_t\dot\psi+\theta\nabla\dot\phi+\nabla\dot\zeta-a\cL\dot\psi=R_\psi,\\
&c_vD_t\dot\zeta+\theta\diver\dot\psi-\kappa a\Delta\dot\zeta=R_\zeta,
\end{aligned}
\right.
\]
with coefficients $(u,\theta,a)$ given by the base solution $U=\bar U+f$, and sources $(R_\phi,R_\psi,R_\zeta)$ satisfying
\begin{equation}\label{eq:general_source_condition}
\norm{R_\phi}_{H^1}+\norm{R_\psi}_2+\norm{R_\zeta}_2
\le C_R(\delta^2+N)\norm V_{H^1}.
\end{equation}
Then for sufficiently small $\delta$ and $N=\norm Z_{H^3}$, the functional $P_1(V)$ in \eqref{eq:Pm} satisfies
\begin{equation}\label{eq:general_physical_bound}
P_1(V)'+c_pD_1(V)
\le C_p(\delta+N)\bigl\{D_1(V)+\norm V_2^2\bigr\},
\end{equation}
\leavevmode{where $c_p,C_p>0$ depend only on the system parameters and $C_R$. By standard Friedrichs mollification, the integrated form of \eqref{eq:general_physical_bound} holds for the time and transverse spatial derivatives of strong solutions.}
\end{lemma}

\begin{proof}
Applying $\partial^\alpha$ ($|\alpha|\le1$) to the system for $V$, testing with $(\theta_*\dot\phi_\alpha,\dot\psi_\alpha,\dot\zeta_\alpha/\theta_*)$, and adding the cross functional $\int\dot\psi\cdot\nabla\dot\phi\dd$ yields the identities \eqref{eq:fixed_symmetric_identity} and \eqref{eq:fixed_compensator_identity} for $V$ with sources $(R_\phi,R_\psi,R_\zeta)$. The transport, pressure, coefficient-gradient, and first-order commutator terms are bounded as in the proof of Lemma~\ref{lem:physical_h1} by
\[
C(\delta+N)\bigl(\norm V_{H^1}^2+D_1(V)\bigr).
\]
For the prescribed sources, integrating by parts via
\begin{gather*}
\sum_{|\alpha|\le1}\ip{\partial^\alpha R_\psi}{\dot\psi_\alpha}=\ip{R_\psi}{\dot\psi-\Delta\dot\psi},\qquad
\sum_{|\alpha|\le1}\ip{\partial^\alpha R_\zeta}{\dot\zeta_\alpha}=\ip{R_\zeta}{\dot\zeta-\Delta\dot\zeta},\\
\int\dot\psi\cdot\nabla R_\phi\dd=-\int(\diver\dot\psi)R_\phi\dd,
\end{gather*}
and applying \eqref{eq:general_source_condition} bounds all source pairings by $C(\delta^2+N)(\norm V_{H^1}^2+D_1(V))$, which yields \eqref{eq:general_physical_bound}.
\end{proof}

\section[Global existence and time decay]{{Global existence and time decay}}\label{sec:global_decay}

\leavevmode{We now combine the negative Sobolev and physical energy estimates to obtain global bounds.}

\begin{proposition}\label{lem:closed_energy}
There exist fixed coupling constants $L_1,L_2,L_3\ge1$ and positive constants $c,C>0$ such that for sufficiently small constants $\delta_1,\varepsilon_1>0$ in \eqref{eq:bootstrap}, the combined energies
\begin{equation}\label{eq:closed}
\begin{gathered}
\cE_m=\cE_-(g)+L_mP_m(Z),\\
\cE_m'+c\bigl\{D_0+D_s+D_m(Z)\bigr\}\le0,
\qquad m=1,2,3,
\end{gathered}
\end{equation}
satisfy the uniform equivalence
\begin{equation}\label{eq:total_equiv}
\cE_m\asymp\norm f_{\Hm}^2+\norm Z_{H^m}^2.
\end{equation}
In particular, on $[0,T]$, it holds that
\begin{equation}\label{eq:apriori}
\begin{split}
\sup_{0\le t\le T}\left\{\norm{f(t)}_{\Hm}^2+\norm{Z(t)}_{H^3}^2\right\}
&+\int_0^T\left\{\norm{Z(t)}_{H^3}^2+\norm{\nabla(\psi,\zeta)(t)}_{H^3}^2\right\}\,dt\\
&\le C\left\{\norm{f_0}_{\Hm}^2+\norm{Z_0}_{H^3}^2\right\},
\end{split}
\end{equation}
with $C>0$ independent of $T$, $\delta$, and $\varepsilon$.
\end{proposition}

\begin{proof}
Inserting the density-gradient recovery bound \eqref{eq:rho_recovery} into \eqref{eq:negative_estimate} and absorbing $C\delta^4(D_0+\delta^2D_s)$ into the left-hand side for small $\delta$, we combine Proposition~\ref{prop:negative_energy}, Lemma~\ref{lem:negative_recovery}, and Proposition~\ref{prop:physical} into
\begin{align}
\cE_-'+c_-(D_0+D_s)
&\le C_-D_m(Z)+C_NN\bigl\{D_0+D_s+D_m(Z)\bigr\},\label{eq:inverse_final}\\
P_m'+c_pD_m(Z)
&\le C_p(\delta+N)\bigl\{D_0+D_s+D_m(Z)\bigr\},\label{eq:physical_final}
\end{align}
where $c_-,c_p,C_-,C_N,C_p>0$ are chosen uniformly across $m=1,2,3$ and the time-derivative estimates in Lemmas~\ref{lem:tangent_physical}--\ref{lem:tangent_inverse}.

Fix $L_m\ge1$ such that $L_mc_p-C_-\ge1$ ($m=1,2,3$), and set $L_*=\max\{L_1,L_2,L_3\}$. Adding $L_m$ times \eqref{eq:physical_final} to \eqref{eq:inverse_final} and choosing $\delta_1,\varepsilon_1>0$ to satisfy
\begin{equation}\label{eq:threshold}
L_*C_p\delta_1+(C_N+L_*C_p)\varepsilon_1
\le\frac12\min\{c_-,1\}
\end{equation}
yields \eqref{eq:closed}.

Since $\norm{\nabla\rr}_{\Hm}=\norm{\rr}_2$ and $g=f+(0,d\nabla\rr,0)^T$, we have
\[
\norm g_{\Hm}\le\norm f_{\Hm}+d\norm{\rr}_2,\qquad
\norm f_{\Hm}\le\norm g_{\Hm}+d\norm{\rr}_2,\qquad
\norm{\rr}_2\le C\norm Z_{H^m},
\]
which together with $\cE_-(g)\asymp\norm g_{\Hm}^2$ and $P_m(Z)\asymp\norm Z_{H^m}^2$ establishes \eqref{eq:total_equiv}. Finally, since
\[
\norm Z_{H^3}^2+\norm{\nabla(\psi,\zeta)}_{H^3}^2\le C\bigl(D_0+D_s+D_3(Z)\bigr)
\]
by \eqref{eq:recovery}, integrating \eqref{eq:closed} at $m=3$ over $[0,T]$ and applying \eqref{eq:total_equiv} proves \eqref{eq:apriori}.
\end{proof}

\begin{proof}[Proof of global existence in Theorem~\ref{thm:main}]
\leavevmode{By the local well-posedness theory recalled in Section~\ref{sec:prelim}, it suffices to close the a priori assumption \eqref{eq:bootstrap} and maintain the uniform positivity of $(\rho,\theta)$. For $\varepsilon_0>0$ sufficiently small, the a priori bound \eqref{eq:apriori} gives $\sup_{t\in[0,T]}\norm{Z(t)}_{H^3}\le C\varepsilon\le\frac12\varepsilon_1$, which by Sobolev embedding also keeps $\rho$ and $\theta$ uniformly bounded above and away from zero. A standard continuity argument then extends the local strong solution to all $t\ge0$, and smooth approximation in $H^3(\R^3)\cap\Hm(\R^3)$ together with continuous dependence and weak lower semicontinuity yields \eqref{eq:global}.}
\end{proof}

\leavevmode{By standard Friedrichs mollification and the evolution equations, the integrated differential inequalities derived below for smooth solutions carry over to the global strong solution.}

To prove the decay estimates \eqref{eq:rate}--\eqref{eq:trate}, we differentiate the perturbation equations with respect to $t$ in the moving frame and write
\begin{equation}\label{eq:tangent_variables}
\dot f=f_t=(\dot\rho,\dot m,\dot E)^T,\qquad
\dot Z=Z_t=(\dot\phi,\dot\psi,\dot\zeta)^T,\qquad
\dot g=g_t=\dot f+(0,d\nabla\dot\rho,0)^T.
\end{equation}

\begin{lemma}\label{lem:tangent_physical}
For sufficiently small $\delta$ and $N=\norm Z_{H^3}$, the time derivative $\dot Z=Z_t$ satisfies
\begin{equation}\label{eq:tangentphysical}
P_1(\dot Z)'+c_pD_1(\dot Z)
\le C_p(\delta+N)\bigl\{D_1(\dot Z)+\norm{\dot Z}_2^2\bigr\}.
\end{equation}
\end{lemma}

\begin{proof}
Differentiating the primitive perturbation system \eqref{eq:physical} with respect to $t$ in the stationary shock frame shows that $\dot Z=(\dot\phi,\dot\psi,\dot\zeta)^T$ solves
\begin{equation}\label{eq:tangent}
\left\{
\begin{aligned}
&D_t\dot\phi+\diver\dot\psi=R_\phi,\\
&D_t\dot\psi+\theta\nabla\dot\phi+\nabla\dot\zeta-a\cL\dot\psi=R_\psi,\\
&c_vD_t\dot\zeta+\theta\diver\dot\psi-\kappa a\Delta\dot\zeta=R_\zeta,
\end{aligned}
\right.
\end{equation}
where
\begin{align}
R_\phi&=-\dot\psi\cdot\nabla\log\rho,\nonumber\\
R_\psi&=-\dot\psi\cdot\nabla u-\dot\zeta\nabla\log\rho-a\dot\phi\cL u,\label{eq:Rs}\\
R_\zeta&=-c_v\dot\psi\cdot\nabla\theta-\dot\zeta\diver u-\kappa a\dot\phi\Delta\theta
-a\dot\phi S(u):\nabla u+2aS(u):\nabla\dot\psi.\nonumber
\end{align}
\leavevmode{Applying $\partial^\alpha$ ($|\alpha|\le1$) to \eqref{eq:tangent} yields a differentiated system of the form \eqref{eq:physical_differentiated} for $(\dot\phi_\alpha,\dot\psi_\alpha,\dot\zeta_\alpha)$, whose right-hand sides $(R_{\phi,\alpha},R_{\psi,\alpha},R_{\zeta,\alpha})$ are given by \eqref{eq:physical_commutators} with $(F_\phi,F_\psi,F_\zeta)$ and $(\phi,\psi,\zeta)$ replaced by $(R_\phi,R_\psi,R_\zeta)$ and $(\dot\phi,\dot\psi,\dot\zeta)$, respectively. Testing the resulting equations as in Lemma~\ref{lem:general_physical} yields the symmetric energy identity}
\begin{align}\label{eq:tangent_symmetric_identity}
S_1(\dot Z)'+\mathcal V_1(\dot Z)
={}&\sum_{|\alpha|\le1}\int\Bigg\{
\frac{\diver u}{2}\left(\theta_*\dot\phi_\alpha^2+|\dot\psi_\alpha|^2+\frac{c_v}{\theta_*}\dot\zeta_\alpha^2\right)
-(\theta-\theta_*)\left(\dot\psi_\alpha\cdot\nabla\dot\phi_\alpha
+\frac{\dot\zeta_\alpha}{\theta_*}\diver\dot\psi_\alpha\right)\nonumber\\
&-\mu\sum_{j=1}^3(\partial_ja)\dot\psi_\alpha\cdot\partial_j\dot\psi_\alpha
-(\mu+\lambda)(\nabla a\cdot\dot\psi_\alpha)\diver\dot\psi_\alpha
-\frac\kappa{\theta_*}\dot\zeta_\alpha\nabla a\cdot\nabla\dot\zeta_\alpha\nonumber\\
&+\theta_*\dot\phi_\alpha R_{\phi,\alpha}
+\dot\psi_\alpha\cdot R_{\psi,\alpha}
+\frac{\dot\zeta_\alpha}{\theta_*}R_{\zeta,\alpha}\Bigg\}\dd,
\end{align}
and the cross energy identity
\begin{align}\label{eq:tangent_compensator_identity}
I_1(\dot Z)'+\int\theta|\nabla\dot\phi|^2\dd
={}&\int(a\cL\dot\psi-\nabla\dot\zeta)\cdot\nabla\dot\phi\dd
+\norm{\diver\dot\psi}_2^2\nonumber\\
&+\int\bigl\{(\diver u)\dot\psi\cdot\nabla\dot\phi
-\dot\psi\cdot(\nabla u)^T\nabla\dot\phi
+R_\psi\cdot\nabla\dot\phi-(\diver\dot\psi)R_\phi\bigr\}\dd.
\end{align}

To apply Lemma~\ref{lem:general_physical}, we verify the source bound \eqref{eq:general_source_condition} for $(R_\phi,R_\psi,R_\zeta)$. Using $\norm{\dot\psi}_6\le C\norm{\dot Z}_{H^1}$ and $\norm{D^2\phi}_3\le CN$, we have
\[
\norm{R_\phi}_{H^1}
\le C(\delta^2+N)\norm{\dot Z}_{H^1}+C\norm{\dot\psi}_6\norm{D^2\phi}_3
\le C(\delta^2+N)\norm{\dot Z}_{H^1}.
\]
In $R_\psi$ and $R_\zeta$, the convective terms are bounded by $C(\delta^2+N)\norm{\dot Z}_{H^1}$, while H\"older's inequality gives
\begin{align*}
\norm{a\dot\phi\cL u}_2&\le C\norm{\dot\phi}_6\norm{\cL\psi}_3+C\norm{\cL\bar u}_\infty\norm{\dot\phi}_2\le CN\norm{\nabla\dot\phi}_2+C\delta^2\norm{\dot\phi}_2,\\
\norm{a\dot\phi\Delta\theta}_2&\le C\norm{\dot\phi}_6\norm{\Delta\zeta}_3+C\norm{\bt''}_\infty\norm{\dot\phi}_2\le CN\norm{\nabla\dot\phi}_2+C\delta^2\norm{\dot\phi}_2.
\end{align*}
Decomposing the viscous dissipation via
\[
S(u):\nabla u=\nu(\bu')^2+2S(\bar u):\nabla\psi+S(\psi):\nabla\psi
\]
similarly yields
\[
\norm{a\dot\phi S(u):\nabla u}_2\le C\delta^4\norm{\dot\phi}_2+C(\delta^2+N)N\norm{\nabla\dot\phi}_2,\qquad
\norm{aS(u):\nabla\dot\psi}_2\le C(\delta^2+N)\norm{\nabla\dot\psi}_2.
\]
Summing these bounds gives
\begin{equation}\label{eq:sourcebound}
\norm{R_\phi}_{H^1}+\norm{R_\psi}_2+\norm{R_\zeta}_2\le C(\delta^2+N)\norm{\dot Z}_{H^1}.
\end{equation}
Combining \eqref{eq:tangent_symmetric_identity} and \eqref{eq:tangent_compensator_identity} yields
\begin{equation}\label{eq:tangentphysical_raw}
\begin{split}
P_1(\dot Z)'+c_pD_1(\dot Z)
\le{}&C(\delta+N)\bigl\{\norm{\dot Z}_{H^1}^2+D_1(\dot Z)\bigr\}\\
&+C\bigl(\norm{R_\phi}_{H^1}+\norm{R_\psi}_2+\norm{R_\zeta}_2\bigr)\bigl\{\norm{\dot Z}_{H^1}+D_1(\dot Z)^{1/2}\bigr\},
\end{split}
\end{equation}
which together with \eqref{eq:sourcebound} and $\norm{\dot Z}_{H^1}^2\le C(D_1(\dot Z)+\norm{\dot Z}_2^2)$ establishes \eqref{eq:tangentphysical}.
\end{proof}

Next, for $\dot g=g_t$ and its potential $\dot h=\cA^{-1}\dot g$, we set
\begin{gather*}
\cE_-(\dot g)=\frac12\ip{\dot g}{\cA^{-1}\dot g},\qquad
D_0(\dot h)=\norm{D^2\dot h}_2^2,\qquad
D_s(\dot h)=\int\omega|\bu'|\,|\nabla\dot h|^2\dd.
\end{gather*}

\begin{lemma}\label{lem:tangent_inverse}
For sufficiently small $\delta$ and $N$, it holds that
\begin{equation}\label{eq:tangent_inverse_bound}
\cE_-(\dot g)'+c_-\bigl\{D_0(\dot h)+D_s(\dot h)\bigr\}
\le C_-D_1(\dot Z)+CN\bigl\{D_0(\dot h)+D_s(\dot h)+D_1(\dot Z)\bigr\}.
\end{equation}
\end{lemma}

\begin{proof}
Differentiating \eqref{eq:effective_perturbation} in $t$ gives
\begin{equation}\label{eq:tangent_effective}
\dot g_t-s\partial_1\dot g+\sum_{j=1}^3\partial_j(A_j\dot g)
=\mathcal B\dot g+\dot\cN_g,
\end{equation}
with
\[
\dot\cN_g=
\begin{pmatrix}
0\\
-\diver\dot{\mathcal Q}_m+\cL\dot{\mathcal Q}_u\\
-\diver\dot{\mathcal Q}_E+\kappa\Delta\dot{\mathcal Q}_\theta
+\diver\{S(\dot{\mathcal Q}_u)\bar u+S(\bar u)\dot{\mathcal Q}_u
+S(\dot\psi)\psi+S(\psi)\dot\psi\}
\end{pmatrix}.
\]
Pairing \eqref{eq:tangent_effective} with $\dot h=\cA^{-1}\dot g$ yields
\begin{equation}\label{eq:tangent_inverse_identity}
\cE_-(\dot g)'=Q_E(\dot h)+V(\dot g,\dot h)+Q_N(\dot g,\dot h),
\end{equation}
where
\[
Q_E(\dot h)=\left\langle\dot g,\sum_{j=1}^3\widetilde A_j^T\partial_j\dot h\right\rangle,\qquad
V(\dot g,\dot h)=\ip{\mathcal B\dot g}{\dot h},\qquad
Q_N(\dot g,\dot h)=\ip{\dot\cN_g}{\dot h}.
\]
Since $\cA$ and $\mathcal B$ are independent of $t$, the Euler and diffusive estimates of Proposition~\ref{prop:negative_energy} apply directly to $(\dot g,\dot h)$, giving
\[
Q_E(\dot h)\le-c_ED_s(\dot h),\qquad
V(\dot g,\dot h)\le-cD_0(\dot h)+\frac{c_E}{2}D_s(\dot h)+C\norm{\nabla\dot\rho}_2^2.
\]

Differentiating the conservative-to-primitive relations in $t$ yields
\begin{equation}\label{eq:change_tangent}
\dot\rho=\rho\dot\phi,\qquad \dot m=\rho(\dot\psi+u\dot\phi),\qquad
\dot E=\rho\left\{c_v\dot\zeta+u\cdot\dot\psi+\left(c_v\theta+\frac12|u|^2\right)\dot\phi\right\}.
\end{equation}
Since $\dot g_\rho=\dot\rho=\rho\dot\phi$ and $\nabla\dot\rho=\rho\nabla\dot\phi+\dot\phi\nabla\rho$, combining \eqref{eq:change_tangent} with the elliptic estimate \eqref{eq:elliptic} for $(\dot g,\dot h)$ gives
\begin{align}
\norm{\dot Z}_2^2+\norm{\dot f}_2^2+\norm{\nabla\dot f}_2^2
&\le C\bigl\{D_0(\dot h)+D_1(\dot Z)+\delta^2D_s(\dot h)\bigr\},\label{eq:tangrec}\\
\norm{\nabla\dot\rho}_2^2
&\le CD_1(\dot Z)+C(N^2+\delta^4)\bigl\{D_0(\dot h)+\delta^2D_s(\dot h)\bigr\}.\nonumber
\end{align}
Differentiating \eqref{eq:remainders} in $t$ shows that
\[
|\dot{\mathcal Q}_u|+|\dot{\mathcal Q}_\theta|+|\dot{\mathcal Q}_F|\le C|Z||\dot Z|,
\]
and hence
\[
\norm{\dot{\mathcal Q}_F}_{6/5}+\norm{\dot{\mathcal Q}_u}_{6/5}\le CN\norm{\dot Z}_2,\qquad
\norm{\dot{\mathcal Q}_u}_2+\norm{\dot{\mathcal Q}_\theta}_2\le CN\norm{\nabla\dot Z}_2.
\]
\leavevmode{For the differentiated viscous work in $\dot\cN_g$, H\"older's inequality gives}
\begingroup\[
\norm{S(\dot\psi)\psi+S(\psi)\dot\psi}_{6/5}
\le C\bigl(\norm\psi_3\norm{\nabla\dot\psi}_2+
\norm{\nabla\psi}_3\norm{\dot\psi}_2\bigr)
\le CN\bigl(\norm{\nabla\dot Z}_2+\norm{\dot Z}_2\bigr).
\]
\endgroup
\leavevmode{Pairing this term with $\nabla\dot h_E\in L^6(\R^3)$ and treating $\diver\{S(\dot{\mathcal Q}_u)\bar u+S(\bar u)\dot{\mathcal Q}_u\}$ via the integration-by-parts identity \eqref{eq:work} with $w=\dot{\mathcal Q}_u$, we deduce from the preceding remainder bounds and \eqref{eq:tangrec} that}
\begin{align}
|Q_N(\dot g,\dot h)|&\le CND_0(\dot h)^{1/2}\bigl(\norm{\dot Z}_2+\norm{\nabla\dot Z}_2\bigr)\nonumber\\
&\le CN\bigl\{D_0(\dot h)+D_1(\dot Z)+\delta^2D_s(\dot h)\bigr\}.
\label{eq:tangentnonlinear}
\end{align}
Inserting these bounds into \eqref{eq:tangent_inverse_identity} and absorbing the $(N^2+\delta^4)$ terms proves \eqref{eq:tangent_inverse_bound}.
\end{proof}

\begin{proposition}\label{lem:tangent_closed}
With $L_1\ge1$ fixed in Proposition~\ref{lem:closed_energy}, the time-derivative energy
\begin{equation}\label{eq:Ew}
{\ETan}:=\cE_-(\dot g)+L_1P_1(\dot Z),\qquad
{\ETan'}+c\bigl\{D_0(\dot h)+D_s(\dot h)+D_1(\dot Z)\bigr\}\le0,
\end{equation}
satisfies
\[
{\ETan}\asymp\norm{\dot f}_{\Hm}^2+\norm{\dot Z}_{H^1}^2.
\]
\end{proposition}

\begin{proof}
Inserting the recovery bound \eqref{eq:tangrec} into \eqref{eq:tangentphysical} and adding $L_1$ times the resulting inequality to \eqref{eq:tangent_inverse_bound} yields \eqref{eq:Ew} under the threshold condition \eqref{eq:threshold}. The norm equivalence follows from
\[
\norm{\dot g-\dot f}_{\Hm}=d\norm{\dot\rho}_2\le C\norm{\dot Z}_{H^1},\qquad
\cE_-(\dot g)\asymp\norm{\dot g}_{\Hm}^2,\qquad
P_1(\dot Z)\asymp\norm{\dot Z}_{H^1}^2.
\]
\end{proof}

\begin{lemma}\label{lem:tangent_decay}
Let $E_0=\norm{f_0}_{\Hm}^2+\norm{Z_0}_{H^3}^2$. Then for all $t\ge0$,
\begin{equation}\label{eq:ewrate}
{\ETan}(t)\le CE_0(1+t)^{-1},
\end{equation}
and
\begin{equation}\label{eq:H1rate}
\norm{Z(t)}_{H^1}^2\le C[-\cE_1'(t)]
\le CE_0(1+t)^{-1/2}.
\end{equation}
Consequently,
\[
\norm{f(t)}_{H^3}+\norm{Z(t)}_{H^3}
\le CE_0^{1/2}(1+t)^{-1/4}.
\]
\end{lemma}

\begin{proof}
From the perturbation equations \eqref{eq:conservative_perturbation} and \eqref{eq:physical}, we have
\begin{equation}\label{eq:timebound}
\norm{Z_t}_{H^1}+\norm{f_t}_{\Hm}\le CN,
\end{equation}
because
\begin{gather*}
\norm{\diver\psi}_{H^1}\le C\norm\psi_{H^2},\qquad
\norm{(a\cL\psi,a\Delta\zeta)}_{H^1}\le C\norm{(\psi,\zeta)}_{H^3},\\
\norm{aS(\psi):\nabla\psi}_{H^1}\le CN^2,\qquad
\norm{f_t}_{\Hm}\le\norm{\mathcal F}_2\le CN.
\end{gather*}
By \eqref{eq:apriori}, it follows that
\[
{\ETan}(0)\le CE_0,\qquad
\int_0^\infty{\ETan}(t)\,dt
\le C\int_0^\infty\norm{Z(t)}_{H^3}^2\,dt\le CE_0.
\]
Since ${\ETan}(t)$ is nonincreasing by \eqref{eq:Ew}, we have
\[
{\ETan}(t)\le{\ETan}(0)\le CE_0,\qquad
t{\ETan}(t)\le\int_0^t{\ETan}(s)\,ds\le CE_0,
\]
which proves \eqref{eq:ewrate}.

Next, consider the stationary positive-definite quadratic form
\[
Q(b,V)=\frac12\ip b{\cA^{-1}b}+L_1P_1(V).
\]
Since $\cE_1(t)=Q(g,Z)$ and ${\ETan}(t)=Q(\dot g,\dot Z)$ have time-independent coefficients, the Cauchy--Schwarz inequality for the bilinear form associated with $Q$ gives
\[
|\cE_1'(t)|\le 2\sqrt{\cE_1(t){\ETan}(t)}.
\]
Combining \eqref{eq:recovery} at $m=1$ with \eqref{eq:closed} and \eqref{eq:ewrate} gives
\[
\norm{Z(t)}_{H^1}^2
\le C\bigl\{D_0+D_s+D_1(Z)\bigr\}
\le-C\cE_1'(t)
\le C\sqrt{\cE_1(t){\ETan}(t)}
\le CE_0(1+t)^{-1/2},
\]
which proves \eqref{eq:H1rate}.

Finally, using $P_3(Z)\asymp\norm Z_{H^3}^2\le C(D_3(Z)+\norm Z_2^2)$ in \eqref{eq:phys_est} yields the damped inequality
\[
P_3(Z)'+\lambda_0P_3(Z)\le C\norm Z_2^2\le CE_0(1+t)^{-1/2}
\]
for a constant $\lambda_0>0$. Integrating this differential inequality and splitting the Duhamel convolution at $t/2$ gives
\[
P_3(Z(t))\le CE_0(1+t)^{-1/2},
\]
\leavevmode{which together with \eqref{eq:positiveequiv} proves the stated decay bounds.}

\leavevmode{Finally, for each fixed global solution, the monotonicity and integrability of $\ETan$ imply}
\begingroup\[
\frac{t}{2}\ETan(t)\le\int_{t/2}^t\ETan(\tau)\,d\tau=o(1)\qquad\text{as }t\to\infty.
\]
\endgroup
\leavevmode{Consequently, $\ETan(t)=o(t^{-1})$ and $\norm{Z(t)}_{H^1}^2=o(t^{-1/2})$ by the preceding interpolation inequality, and splitting the Duhamel integral for $P_3(Z(t))$ at $t/2$ yields $P_3(Z(t))=o(t^{-1/2})$, which establishes the asymptotic relations in Remark~\ref{rem:data}.}
\end{proof}

\begin{corollary}\label{cor:transverse}
Under the hypotheses of Theorem~\ref{thm:main}, the transverse spatial derivatives satisfy
\[
\norm{\nabla_\perp Z(t)}_{H^1}+\norm{\nabla_\perp f(t)}_{\Hm}
\le C\varepsilon(1+t)^{-1/2}\qquad\text{for all }t\ge0.
\]
\end{corollary}

\begin{proof}
\leavevmode{For $j\in\{2,3\}$, define the transverse derivative energy}
\begingroup\[
\widetilde\cE_{1,j}(t)
=\cE_-(\partial_jg)+L_1P_1(\partial_jZ)
\asymp\norm{\partial_jf}_{\Hm}^2+\norm{\partial_jZ}_{H^1}^2.
\]
\endgroup
\leavevmode{Since the planar profile $\bar U(x_1)$ is independent of $x_j$ and $[\partial_j,\cA]=0$, differentiating the perturbation system with respect to $x_j$ and repeating the proof of Proposition~\ref{lem:tangent_closed} shows that $\widetilde\cE_{1,j}(t)$ is nonincreasing in $t$. Moreover, since $\norm{\partial_jf}_{\Hm}\le\norm f_2$, the global bound \eqref{eq:global} gives}
\begingroup\[
\widetilde\cE_{1,j}(0)
+\int_0^\infty\widetilde\cE_{1,j}(t)\,dt\le CE_0,
\]
\endgroup
{whence $\widetilde\cE_{1,j}(t)\le CE_0(1+t)^{-1}$ by the same monotonicity argument as in \eqref{eq:ewrate}.}
\end{proof}

\section*{Declaration of AI Assistance}
The central idea of studying the stability of planar shock waves in the homogeneous Sobolev space $\Hm(\R^3)$ originated with the authors. ChatGPT (OpenAI) and Gemini (Google) provided auxiliary assistance during the development of this work. The authors led the research and take full responsibility for all mathematical derivations, results, and conclusions presented in this paper.

\end{document}